\documentclass[openany, a4paper, 11pt]{article}
\usepackage{mathrsfs}
\usepackage{amsfonts}

 \usepackage{mathrsfs,amsfonts,amsmath}
 \usepackage{color}
 \newtheorem{definition}{Definition}[section]
 \newtheorem{hypothesis}{Hypothesis}[section]
 \newtheorem{lemma}{Lemma}[section]
 \newtheorem{proposition}{Proposition}[section]
 \newtheorem{theorem}{Theorem}[section]
 \newtheorem{corollary}{Corollary}[section]
 \newtheorem{Remark}{Remark}[section]
 \newtheorem{example}{Example}[section]

 \def\blemma{\begin{lemma}\sl{}\def\elemma{\end{lemma}}}
 
 \def\btheorem{\begin{theorem}\sl{}\def\etheorem{\end{theorem}}}
 \def\bcorollary{\begin{corollary}\sl{}\def\ecorollary{\end{corollary}}}
 \def\bremark{\begin{Remark}\rm{}\def\eremark{\end{Remark}}}
 \def\bexample{\begin{example}\rm{}\def\eexample{\end{example}}}

 \def\beqlb{\begin{eqnarray}}\def\eeqlb{\end{eqnarray}}
 \def\beqnn{\begin{eqnarray*}}\def\eeqnn{\end{eqnarray*}}

 \def\ar{\!\!\!&}

 \def\mbb{\mathbb}

 \def\proof{\noindent{\it Proof.~~}}\def\qed{\hfill$\Box$\medskip}

 \def\bfE{{\mbox{\bf E}}}\def\bfP{{\mbox{\bf P}}}
 \def\bfvar{{\mbox{\bf var}}}

\begin{document}

\noindent{\small 09fluclim (2009/08)}

\vskip1.0cm

\centerline{\Large\textbf{A Fluctuation Limit Theorem of}}

\smallskip

\centerline{\Large\textbf{Branching Processes with Immigration}}

\smallskip

\centerline{\Large\textbf{and Statistical Applications}}

\medskip\bigskip

\centerline{Chunhua Ma}

\centerline{School of Mathematical Sciences and LPMC,}

\centerline{Nankai University,}

\centerline{Tianjin 300071, P.\,R.\ China}

\centerline{E-mail: {\tt mach@nankai.edu.cn}}

\medskip

\bigskip

\noindent\textit{Abstract.} We prove a general fluctuation limit
theorem for Galton-Watson branching processes with immigration. The
limit is a time-inhomogeneous OU type process driven by a spectrally
positive L\'{e}vy process. As applications of this result, we obtain
some asymptotic estimates for the conditional least squares
estimators of the means and variances of the offspring and
immigration distributions.

\medskip

\noindent\textit{Key words.} Branching process with immigration;
Ornstein-Uhlenbeck type process; one-sided stable process;
fluctuation limit; Poisson random measure;  conditional least
squares estimator

\medskip

\noindent\textit{AMS 2000 subject classifications.} Primary 60J35;
secondary 60J80, 60H20, 60K37.


\bigskip\medskip

\section{Introduction}

\setcounter{equation}{0}

Consider the \textit{Galton-Watson branching process with
immigration}, $\{y(k): k=1,2,\cdots\}$, defined by
 \beqlb
 y(k) = \sum^{y(k-1)}_{j=1}\xi(k,j) + \eta(k),\quad k\geq1,\quad\quad
 y(0)=0,
 \label{GWI-equation}
 \eeqlb
where $\{\xi(k,j): k,j=1,2,\cdots\}$ and $\{\eta(k): k=1,2,\cdots\}$
are two independent families of i.i.d.\ random variables taking
values in $\mathbb{N} :=\{0,1,2,\cdots\}$. The distribution of
$\xi(k,j)$ is called the offspring distribution and the distribution
of $\eta(k)$ is called the immigration distribution. Let $g(\cdot)$
and $h(\cdot)$ be the generating functions of $\xi(k,j)$ and
$\eta(k)$, respectively. It is easy to see that $\{y(k)\}$ is a
discrete-time Markov chain with values in $\mathbb{N}$ and one-step
transition matrix $P(i,j)$ given by
 \beqlb
 \sum_{j=0}^{\infty}P(i,j)s^j=g(s)^ih(s),\qquad i\in\mathbb{N},\
 0\leq s\leq 1.
 \label{trans matrix}
 \eeqlb
For simplicity, we also call $\{y(k)\}$ a \textit{GWI-process with
parameters} $(g,h)$. Assume that the offspring mean $m:=g'(1)$ is
finite. The cases $m>1$, $m=1$ and $m<1$ are referred to
respectively as \textit{supercritical, critical, and subcritical}. A
sequence of GWI-processes $\{y_n(\cdot)\}$ with $(g_n,h_n)$ is said
to be \textit{nearly critical} if $m_n:=g_n'(1)$ converges to $1$ as
$n$ tends to $\infty$.

The estimation problem for the offspring and immigration parameters
in the GWI-process has been extensively studied; see Heyde and
Seneta \cite{HS72, HS74}, Wei and Winnicki \cite{WW90} and the
references therein. It is well known that the conditional least
squares estimators (CLSE), first obtained by Klimko and Nelson
\cite{KN78}, can be used to estimate the offspring mean on the basis
of the observing information on $\{y(k)\}$; see also \cite{HS72} and
\cite{WW90} for other closely related estimators. In the
non-critical case, the CLSE of the offspring mean is consistent and
asymptotically normal (see \cite{KN78}, \cite{WW89,WW90}). However,
it was shown by \cite{S94} and \cite{WW89} that in the critical or
nearly critical case the CLSE is not asymptotically normal. In fact,
when the process is nearly critical and the offspring variance tends
to a positive real number,  Sriram \cite{S94} gave the weak
convergence of GWI-processes to the branching diffusion with
immigration. As a result, the above CLSE of the offspring mean has
the asymptotic distribution which is expressed in terms of the limit
process and the normalizing factor is $n$. Motivated by the similar
statistical application, Isp\'{a}ny \textit{et al.}\;\cite{I05} have
recently obtained a fluctuation limit theorem for the nearly
critical processes where the offspring variances tend to $0$. Such
limit is a time-inhomogeneous Ornstein-Uhlenbeck (OU) processes
driven by a Wiener process. As a consequence, they proved the
asymptotic normality of CLSE of the offspring mean with normalizing
factor $n^{3/2}$. Obviously, the asymptotic behavior of the CLSE in
the critical or nearly critical case is closely related to the limit
theorems of the GWI-processes.

The main objective of this paper is to give a general fluctuation
limit theorem and its applications for processes that allow the
offspring and immigration distributions to have infinite variances.
Fluctuation limits for branching models with immigration have been
investigated by Dawson and Li \cite{DL06}, Isp\'{a}ny \textit{et
al.}\;\cite{I05}, Li \cite{L00} and Li and Ma \cite{LM08}; see also
Dawson \textit{et al.}\;\cite{DFG89} for the type of limits in the
measure-valued setting. In the present paper, we shall consider a
sequence of nearly critical GWI-processes $\{y_n(\cdot)\}$ with
$(g_n,h_n)$ satisfying a set of conditions similar to that of
\cite{L00}. Let us define the sequence
 $
Y_n(t)=y_n([nt])
 $
and consider the rescaled centralized process $Z_n(\cdot)=c_n^{-1}
(Y_n(\cdot)-E[Y_n(\cdot)])$ with certain sequence of positive
constants $c_n$. It turns out that $Z_n(\cdot)$ converges to a
time-inhomogeneous OU type process driven by a spectrally positive
L\'{e}vy process (Theorems \ref{OU limit}). Based on this
fluctuation limit, we show that non-degenerate limit laws still
exist for the above CLSE estimates of means (Theorem
\ref{asymptotic}). Of special interest is the case when the
offspring and immigration distributions belong to the domain of
attraction of a stable law
 with exponent $\alpha$ $(1<\alpha\leq2)$. For simplicity,
suppose that
 \beqnn
 g_n(s)= s+\frac{\gamma}{n}(1-s)^{\alpha} \quad\mbox{and}\quad
  h_n(s) = s+\varpi(1-s)^{\alpha},
  \label{intro-example}
 \eeqnn
where $0<\gamma, \varpi\leq1/\alpha$. Note that for $1<\alpha<2$,
$g_n$ has infinite variance but its heavy-tailed effect weakens as
$n\rightarrow\infty$; for $\alpha=2$, the offspring variance is
$2\gamma/n$ and tends to 0. Then $Z_n(\cdot)$ with
$c_n=n^{1/\alpha}$ converges to a OU type process driven by a
$\alpha$-stable process (Corollary~\ref{example2}). As a
consequence, the CLSE of the offspring mean is asymptotic to a
$\alpha$-stable distribution and the normalizing factor is
$n^{\frac{2\alpha-1}{\alpha}}$ (Corollary \ref{example3}). As
mentioned above, the estimation for the offspring mean in
GWI-process have been systematically studied by \cite{WW89,WW90},
\cite{S94} and \cite{I05}, provided that the offspring variances are
finite. Our results can be regarded as an attempt in the case when
the above assumption fail to hold.

Another interesting case, related to our limit theorem, is that the
offspring variances are finite and tend to $0$, but the offspring
distributions do not satisfy the Lindeberg conditions required in
\cite{I05}. Then the resulting fluctuation limit $Z(\cdot)$ is a OU
type process with positive jumps instead of OU diffusion (Corollary
\ref{example1}). In this case, it is also possible to consider the
CLSE estimates for the offspring and immigration variances. We show
that the CLSE of the offspring variance is consistent and its
asymptotic distribution (with normalizing factor $n$) is expressed
in terms of the jumps of $Z(\cdot)$, while the CLSE of the
immigration variance is not consistent (Theorme \ref{asymptotic
variance}). However, if we return to the case of \cite{I05} (see
also Example \ref{example of diffusion}), the above asymptotic
distribution is degenerate to $0$ and the above immigration variance
estimator becomes consistent (Remark \ref{remark 3.1}). Hence, in
this case, by adding certain conditions on fourth moments we further
prove that these estimators of the offspring and immigration
variances are asymptotically normal with the normalizing factors
$n^{3/2}$ and $n^{1/2}$, respectively (Theorem \ref{diffusion
variance theorem}). This result also contrasts with the
critical-mean and positive-variance case of Winnicki \cite{W91}, in
which the CLSE of the offspring variance is not asymptotically
normal, although it has another limit law with the normalizing
factor $n^{1/2}$.

The remainder of this paper is organized as follows. The main limit
theorems and some examples will be given in section 2. In Section 3
we obtain some asymptotic estimates for the statistics of the
GWI-process, as applications of our limit theorems. Section 4 is
devoted to the proofs of Theorem \ref{Theorem 1}-\ref{OU limit} and
Theorem
\ref{asymptotic2.3}-\ref{diffusion variance theorem}.\\

\noindent\textbf{Notation.} Let $\mbb{R}_+ = [0, \infty)$. For $x\in
\mathbb{R}$, set $\chi(x)= (1 \wedge x ) \vee (-1)$.
$``\overset{p}{\longrightarrow}" $ and
$``\overset{d}{\longrightarrow}" $ denote the convergence of random
variables in probability and convergence in distribution,
respectively. We also make the convention that
 $
\int_r^t = - \int_t^r = \int_{(r,t]}
 $
 and
 $
\int_r^\infty = \int_{(r,\infty)}
 $
for $r\le t \in \mbb{R}$.


\section{Limit theorems and examples}

\setcounter{equation}{0}

Let us consider a sequence of GWI-processes $y_n(\cdot)$ with
parameters $(g_n, h_n)$. A realization of $y_n(\cdot)$ is defined by
\beqlb
 y_n(k) = \sum^{y_n(k-1)}_{j=1}\xi_n(k,j) + \eta_n(k),\quad k\geq1,\quad\quad
 y_n(0)=0,
 \label{GWI-equationn}
 \eeqlb
where $\{\xi_n(k,j)\}$ and $\{\eta_n(k)\}$ are given as in
(\ref{GWI-equation}), but depend on the index $n$. Also, $g_n$ and
$h_n$ are the generating functions of $\xi_n(k,j)$ and $\eta_n(k)$.
Now introduce the sequence
 \beqnn
Y_n(t):=y_n([nt]),\qquad t\geq0,
 \label{y_n}
 \eeqnn
where $[nt]$ denotes the integer-part of $nt$, and
$Y'_n(t):=\sum^{[nt]}_{k=1}\eta_n(k)$. We first prove a limit
theorem for the sequence $(Y_n(\cdot), Y_n'(\cdot))$. Such theorem
is the modification of Theorem 2.1 in \cite{L05}. Let $\{b_n\}$ be a
sequence of positive numbers such that $b_n\rightarrow\infty$ as
$n\rightarrow\infty$. For $0\leq\lambda\leq b_n$, set \beqlb
 R_n(\lambda)=nb_n[(1-\lambda/b_n)-g_n(1-\lambda/b_n)]
 \label{R_n}
 \eeqlb
and
 \beqlb
 F_n(\lambda)=n[1-h_n(1-\lambda/b_n)].
 \label{F_n}
 \eeqlb
Consider the following set of conditions:
\begin{enumerate}
\item[(A)]\;
The sequence $\{R_n\}$ is uniformly Lipschitz on each bounded
interval and converges to a continuous function as $n \rightarrow
\infty$;

\item[(B)]\;
The sequence $\{F_n\}$ converges to a continuous function as $n
\rightarrow \infty$.
 \end{enumerate}

\blemma\label{prop2.1}\; Under condition (A), the limit function $R$
of $\{R_n\}$ has representation
 \beqlb
R(\lambda)= c\lambda -\theta\lambda^2 - \int_0^{\infty} (e^{-\lambda
u} -1+ \lambda u)\, \Lambda_1(du),
 \label{R}
  \eeqlb
where $c\in\mathbb{R}$, $\theta\geq0$, and $\Lambda_1(du)$ is a
$\sigma$-finite measure on $(0,\infty)$ with $\int_0^\infty (u\wedge
u^2)\Lambda_1(du)<\infty$.
 \elemma

\blemma\label{prop2.2}\;Under condition (B), the limit function $F$
of $\{F_n\}$ has representation
 \beqlb
F(\lambda) = d\lambda + \int_0^\infty (1-e^{-\lambda u})
\,\Lambda_2(du),
  \label{F}
  \eeqlb
where $d\geq0$, and $\Lambda_2(dz)$ is a $\sigma$-finite measure on
$(0,\infty)$ with $\int_0^\infty (1\wedge u)\,\Lambda_2(du)<\infty$.
 \elemma
 \btheorem\label{Theorem 1}\; Suppose that (A) and (B) are satisfied.
Then $\big(Y_n(\cdot)\big/b_n, Y'_n(\cdot)\big/b_n\big)$ converges
in distribution on $D([0,\infty), \mathbb{R}_+^2)$ to a
two-dimensional non-negative Markov process $\big(Y(\cdot),
Y'(\cdot)\big)$ with initial value $(0,0)$ and transition semigroup
$(P_t)_{t\geq0}$ given by
 \beqlb
\int_{\mathbb{R}_+^2}e^{-\langle z,u\rangle
}P_t(x,du)=\exp\bigg\{-x_1\psi_t(z_1)-x_2z_2-\int_0^t
F(\psi_s(z_1)+z_2)ds\bigg\},
 \label{semigroup}
 \eeqlb
where $x=(x_1,x_2)\in\mathbb{R}_+^2$, $z=(z_1,z_2)\in\mathbb{R}_+^2$
and $\psi_t(z_1)$ is the unique solution of
 \beqlb
 \frac{d\psi_t}{dt}(z_1)=R(\psi_t(z_1)),
 \qquad\psi_0(z_1)=z_1.
 \label{Ric}
 \eeqlb
  \etheorem
 \bremark\label{remark1}\; {\rm (i)} The process $Y(\cdot)$ is a conservative
\textit{continuous state branching process with immigration}
(CBI-process) and $Y'(\cdot)$ is the immigration part of $Y(\cdot)$.
See Kawazu and Watanabe \cite{KW71} for a complete characterization
of the class of CBI-processes. Furthermore, $(Y(\cdot),Y'(\cdot))$
is a special case of two-dimensional CBI-processes; see \cite{M09}
and the references therein.

{\rm (ii)} Lemma \ref{prop2.1}-\ref{prop2.4} are closely related to
the L\'{e}vy-Khinchin type representations of some class of
continuous functions. Li \cite{L91} used a simple method based on
Bernstein polynomials to prove Lemma \ref{prop2.1} and
\ref{prop2.2}. Here, inspired by Venttsel' \cite{V59}, these lemmas
can be proved in a different but more intuitive way. From the proof
of Lemma \ref{lemma2.3}, we can see that either (A) or (D2) implies
the convergence of the sums of the triangular array of i.i.d.
variables $\{\frac{\xi_n(k,j)-1}{b_n}: k,j=1,2\cdots\}$ or
$\{\frac{\xi_n(k,j)-1}{c_n}: k,j=1,2\cdots\}$. See Grimvall
\cite{G74} for a similar consideration.
 \eremark
\bcorollary\label{example(1)}\,(\cite{KW71}) Let $L(x)$ and $L^*(x)$
be positive functions slowly varying at $\infty$ such that $L(x)\sim
L^*(x)$ as $x\rightarrow\infty$. Consider a sequence of
GWI-processes with
 $(g_n, h_n)$ given by
  \beqnn
  g_n(s) \ar\equiv\ar
  s+\gamma(1-s)^{\alpha}L\Big(\frac{1}{1-s}\Big),\\
  h_n(s) \ar\equiv\ar
  1-\varpi(1-s)^{\alpha-1}L^*\Big(\frac{1}{1-s}\Big),
  \eeqnn
where $1<\alpha\leq2$, $\gamma>0$, $\varpi>0$. Let $b_n$ be the
sequence satisfying
 \beqlb
 b_n\sim [nL(b_n)]^{1/(\alpha-1)}\quad (\sim[nL^*(b_n)]^{1/(\alpha-1)}).
 \label{b_nL}
 \eeqlb
Then $(Y_n(\cdot)/b_n, Y'_n(\cdot)/b_n)$ converges in distribution
on $D([0,\infty), \mathbb{R}_+^2)$ to $(Y(\cdot),Y'(\cdot))$ defined
by
 \beqlb
 dY(t)=\sqrt[\alpha]{Y(t-)}\,dX(t)+dY'(t), \qquad
 Y(0)=Y'(0)=0,
 \label{stable1}
 \eeqlb
where $X(t)$ is a spectrally positive $\alpha$-\,stable L\'{e}vy
process with Laplace exponent
 $
R(\lambda)=-\gamma\lambda^{\alpha},
 $
and $Y'(t)$ is a $(\alpha-1)$-\,stable subordinator with Laplace
exponent
 $F(\lambda)=\varpi\lambda^{\alpha-1}$, independent of $X$.
\ecorollary

\bremark Sometimes the above process can be regarded as the
branching process conditioned on not being extinct in the distant
future, or $Q$-process; see Lambert \cite{L07}. The pathwise
uniqueness for the type of SDE (\ref{stable1}) has recently been
proved by Fu and Li \cite{FL08}. \eremark

{\it Proof of Corollary \ref{example(1)}\;} Without loss of
generality, consider $R_n(\lambda)$ by (\ref{R_n})  on
$\lambda\in[0,1]$. It is easy to see that $
\lim_{n\rightarrow\infty}R_n(\lambda)=-\gamma\lambda^{\alpha}.
 $
For $\lambda>0$,
$|R'_n(\lambda)|=\alpha\gamma\lambda^{\alpha-1}[nL^{\diamond}(b_n/\lambda)/b_n^{\alpha-1}]$,
where $L^{\diamond}(x)$ is a slowly varying functions such that
$L^\diamond(x)\sim L(x)$, as $x\rightarrow\infty$ (cf. \cite[Theorem
1.8.2]{BGT87}). By the representation theorem of $L^\diamond(\cdot)$
(\cite[Theorem 1.3.1]{BGT87}),
 \beqnn
L^\diamond(b_n/\lambda)/L^\diamond(b_n)=\{q(b_n/\lambda)/q(b_n)\}\exp
\bigg\{\int^{b_n/\lambda}_{b_n}\epsilon(u)du/u\bigg\},
 \eeqnn
where $q(x)\rightarrow q\in(0,\infty)$, $\epsilon(x)\rightarrow 0$
as $x\rightarrow\infty$. Fix $0<\varepsilon<q\wedge(\alpha-1)$.
There exists $x_0>0$ such that $q-\varepsilon<q(x)<q+\varepsilon$
and $|\epsilon(x)|<\varepsilon$, if $x>x_0$. Note that
$b_n/\lambda\geq b_n$ for $0<\lambda<1$ and choose sufficiently
large $n$, we have
 \beqnn
 \lambda^{\alpha-1}[L^\diamond(b_n/\lambda)/L^\diamond(b_n)]\leq
 \frac{q+\varepsilon}{q-\varepsilon}\lambda^{\alpha-1-\varepsilon}.
 \eeqnn
Then, by the above inequality and (\ref{b_nL}),
$\sup_n|R'_n(\lambda)|$ is bounded in $\lambda\in(0,1]$. Note that
$R'_n(0)=0$ and thus (A) holds. Also, it is not hard to see that
 $
 \lim_{n\rightarrow\infty}F_n(\lambda)=\varpi\lambda^{\alpha-1}.
 $
By Theorem \ref{Theorem 1}, the limit process $(Y(\cdot),Y'(\cdot))$
is defined by (\ref{semigroup}) and (\ref{Ric}) with
$R(\lambda)=-\gamma\lambda^\alpha$ and
$F(\lambda)=\varpi\lambda^{\alpha-1}$. By \cite{FL08},
$(Y(\cdot),Y'(\cdot))$ is the unique solution of the above
stochastic equation system.\qed

Now we turn to study the fluctuation limit for the sequence
$Y_n(\cdot)$. Assume that $m_n=g'_n(1)$ and $\omega_n=h'_n(1)$ are
finite. Let $\{c_n\}$ be a sequence of positive numbers. Set
 \beqlb
 G_n(\lambda)=n^2[(1-m_n\lambda/c_n)-g_n(1-\lambda/c_n)]
 \label{G_n}
 \eeqlb
and
 \beqlb
 H_n(\lambda)=n[(1-\omega_n\lambda/c_n)-h_n(1-\lambda/c_n)],
 \label{H_n}
 \eeqlb
for $0\leq\lambda\leq c_n$. We will need the following conditions:
 \begin{enumerate}

\item[(C)]\;
$n/c_n\rightarrow\infty$ and $n/c_n^2\rightarrow\gamma_{0}$ as
$n\rightarrow\infty$, for some $\gamma_0\geq0$;

\item[(D1)]\;
$n(m_n-1)\rightarrow a$ as $n\rightarrow\infty$, for some
$a\in\mathbb{R}$;

\item[(D2)]\;
The sequence $\{G_n\}$ is uniformly Lipschitz on each bounded
interval and converges to a continuous function as $n \rightarrow
\infty$;

\item[(E1)]\;
$\omega_n\rightarrow \omega$ as $n\rightarrow\infty$, for some
$\omega\geq0$;

\item[(E2)]\;
The sequence $\{H_n\}$ is uniformly Lipschitz on each bounded
interval and converges to a continuous function as $n \rightarrow
\infty$.
 \end{enumerate}
\blemma\label{lemma2.3}\; Under conditions (C) and (D1,2), the limit
function $G$ of $\{G_n\}$ has representation
 \beqlb
G(\lambda)= \beta_1\lambda -\sigma_1\lambda^2 - \int_0^{\infty}
(e^{-\lambda u} -1+ \lambda u)\mu(du),
 \label{G}
  \eeqlb
where $\beta_1\in\mathbb{R}$, $\sigma_1\geq0$ and $2\sigma_1\geq
a\gamma_0$, and $\mu(du)$ is a $\sigma$-finite measure on
$(0,\infty)$ with $\int_0^\infty (u\wedge u^2)\mu(du)<\infty$.
 \elemma

\blemma\label{prop2.4}\;Under conditions (C) and (E1,2), the limit
function $H$ of $\{H_n\}$ has representation
 \beqlb
H(\lambda) = \beta_2\lambda -\sigma_2\lambda^2- \int_0^\infty
(e^{-\lambda u} -1+ \lambda u)\nu(du),
  \label{H}
  \eeqlb
where $\beta_2\in\mathbb{R}$, $\sigma_2\geq0$ and
$2\sigma_2+\omega\gamma_0\geq\omega^2\gamma_0$, and $\nu(du)$ is a
$\sigma$-finite measure on $(0,\infty)$ with $\int_0^\infty (u\wedge
u^2)\nu(du)<\infty$.
 \elemma
Let $\phi(s)=\omega\int_0^se^{au}du$ and let
$\varrho(s)=(2\sigma_1-a\gamma_0)\phi(s)+2\sigma_2+\omega(1-\omega)\gamma_0$.
By the above representations, $\varrho(s)\geq0$ for $s\geq0$. We
actually obtain a set of parameters
$(\beta_1,\beta_2,\varrho(\cdot),\phi(\cdot),\mu,\nu)$ which will be
used to characterizes our limit processes. Let $Z_n(\cdot)$ be
defined by
 \beqlb
Z_n(t)=\frac{Y_n(t)-\bfE[Y_n(t)]}{c_n}.
 \label{Z_n}
  \eeqlb
Our main result of this paper is the following fluctuation limit
theorem.

\btheorem\label{OU limit}\; Suppose that conditions (C), (D1,2) and
(E1,2) are satisfied. Let $B(t)$ be a one-dimensional Brownian
motion, $N_0(ds,du)$ be a Poisson random measure on $(0,\infty)
\times \mathbb{R}_+$ with intensity $ds\nu(du)$ and
$N_1(ds,du,d\zeta)$ be a Poisson random measure on $(0,\infty)
\times \mathbb{R}_+\times(0,\infty)$ with intensity $ds\mu(du)
d\zeta$. Suppose that $B$, $N_0$ and $N_1$ are independent of each
other. Then $Z_n(\cdot)$ converges in distribution on $D([0,\infty),
\mathbb{R})$ to a time-inhomogeneous OU type process $Z(\cdot)$,
which can be constructed as the unique solution of the following
stochastic equation
 \beqlb\label{OU}
Z(t)
 \ar=\ar
 \int_0^t\big(\beta_2 + \beta_1\phi(s) + aZ(s)\big)\,ds +
\int_0^t\sqrt{\varrho(s)}\,dB(s)\nonumber\\
 \ar \ar
+\int_0^t\int_{\mathbb{R}_+}u\tilde{N}_0(ds,du)
+\int_0^t\int_{\mathbb{R}_+}\int_0^{\phi(s)}u
\tilde{N}_1(ds,du,d\zeta),
 \label{OU equation}
 \eeqlb
where $\tilde{N}_0(ds,du)=N_0(ds,du)-ds\nu(du)$ and
$\tilde{N}_1(ds,du,d\zeta)= N_1(ds,du,d\zeta)-ds\mu(du)d\zeta$.
 \etheorem

\bremark\label{remark 2.1}\; The conditions of Theorem \ref{OU
limit} imply that $Y_n(t)/n$ converges weakly to the deterministic
function $\phi(t)=\omega\int_0^te^{as}ds$. In fact, consider $R_n$
in (\ref{R_n}) and $F_n$ in (\ref{F_n}) with $b_n=n$. Note that
$R_n(\lambda)=G_n(c_n\lambda/n)+n(m_n-1)\lambda$. Then by conditions
(C), (D1,2) and Lemma \ref{prop2.1}, $R'_n(\lambda)$ is uniformly
bounded in each bounced interval and
$\lim_{n\rightarrow\infty}R_n(\lambda)=a\lambda$. In a similar way,
we also have $\lim_{n\rightarrow\infty}F_n(\lambda)=\omega\lambda$.
The above weak convergence result follows from Theorem \ref{Theorem
1}.
 \eremark

\bcorollary\label{example1}\, Let $\{y_n(k)\}$ be defined as in the
beginning of this section. In addition to conditions (D1) and (E1),
we assume that $\pi_n=\bfvar\,\xi_n(1,1)< \infty$,
$r_n=\bfvar\,\eta_n(1)<\infty$, and the following conditions hold:
 \begin{itemize}

 \item[(a.1)]\;
The sequence
$\tilde{\mu}_n(\cdot)=n\bfE\Big[(\xi_n(1,1)-m_n)^2{\bf1}_{\{(\xi_n(1,1)-m_n)/\sqrt{n}\
\in\,\cdot\,\}}\Big]$ converges weakly to a finite measure denoted
by $\tilde{\mu}(\cdot)$, as $n\rightarrow\infty$;

 \item[(a.2)]\;
The sequence
$\tilde{\nu}_n(\cdot)=\bfE\Big[(\eta_n(1)-\omega_n)^2{\bf
1}_{\{(\eta_n(1)-\omega_n)/\sqrt{n}\ \in\,\cdot\,\}}\Big]$ converges
weakly to a finite measure denoted by $\tilde{\nu}(\cdot)$, as
$n\rightarrow\infty$.

 \end{itemize}
Let $Z_n(\cdot)$ be defined by (\ref{Z_n}) with $c_n=\sqrt{n}$. Then
$Z_n(\cdot)$ converges in distribution on $D([0,\infty),
\mathbb{R})$ to a OU type process $Z(\cdot)$ whose L\'{e}vy measure
has finite second moment, i.e.
 \beqlb\label{OU_square}
Z(t)
 \ar=\ar
a\int_0^tZ(s)\,ds + \int_0^t\sqrt{\varrho(s)}\,dB(s)
+\int_0^t\int_{\mathbb{R}_+}u\tilde{N}_0(ds,du)\nonumber\\
\ar \ar +\int_0^t\int_{\mathbb{R}_+}\int_0^{\phi(s)}u
\tilde{N}_1(ds,du,d\zeta),
 \label{OU_square equation}
 \eeqlb
where $\varrho(s)=\tilde{\nu}(\{0\})+\tilde{\mu}(\{0\})\phi(s)$ and
similarly $\phi(s)=\omega\int_0^se^{au}du$. $B$, $N_0$ and $N_1$ are
defined as in (\ref{OU equation}), but the corresponding intensities
of $N_0$ and $N_1$ are given by $ds\nu(du)$ and $ds\mu(du)d\zeta$
with $\nu(du)=u^{-2}{\bf1}_{\{u>0\}}\tilde{\nu}(du)$ and
$\mu(du)=u^{-2}{\bf1}_{\{u>0\}}\tilde{\mu}(du)$.
 \ecorollary

\proof  By Theorem \ref{OU limit}, it suffices to check (D2) and
(E2) are satisfied. First it follows from (a.1,2) that
$\tilde{\mu}(\cdot)$ and $\tilde{\nu}(\cdot)$ are supported by
$[0,\infty)$. Consider $G_n(\lambda)$ in (\ref{G_n}) with
$c_n=\sqrt{n}$. Without loss of generality, we restrict ourselves to
$\lambda\in[0,1]$. Making a Taylor expansion of $g_n$ about $1$, we
have
$G_n(\lambda)=-n\int_0^1(1-s)g''_n(1-s\lambda/\sqrt{n}\,)\lambda^2ds$.
Note that
 \beqnn
 ng_n''(1-s\lambda/\sqrt{n})\ar=\ar
  (1-s\lambda/\sqrt{n})^{m_n-2}\Big(\int e^{u\sqrt{n}\ln^{(1-s\lambda/\sqrt{n})}}\tilde{\mu}_n(du)
  \nonumber\\
  \ar\ar+\,nm_n(m_n-1)\bfE[(1-s\lambda/\sqrt{n})^{\xi_n(1,1)-m_n}]\nonumber\\
  \ar\ar+\,n(2m_n-1)\bfE[(\xi_n(1,1)-m_n)(1-s\lambda/\sqrt{n})^{\xi_n(1,1)-m_n}]\Big).
 \eeqnn
Fix $s,\lambda\in(0,1]$ and choose sufficiently large $n$ such that
$0<s\lambda/\sqrt{n}\leq1/2$. We have that
 \beqnn
 \Big|\int e^{u\sqrt{n}\ln^{(1-s\lambda/\sqrt{n})}}\tilde{\mu}_n(du)-
 \int e^{-us\lambda}\tilde{\mu}_n(du)\Big|\leq
 \int (e^{-us\lambda}|u|/\sqrt{n})\vee(e(e^{1/\sqrt{n}\,}-1))\tilde{\mu}_n(du),
 \eeqnn
 \beqnn
 \big|1-(1-s\lambda/\sqrt{n})^{\xi_n(1,1)-m_n}\big|\leq2e^2
 |\xi_n(1,1)-m_n|/\sqrt{n}.
 \eeqnn
Note that $\tilde{\mu}_n(\cdot)$ is supported by
$\{(k-m_n)/\sqrt{n}:k=0,1,\cdots\}$. By conditions (D1), (E1) and
(a.1), it is easy to see that $ng''(1-s\lambda/\sqrt{n})\rightarrow
a+\int_{[0,\infty)}e^{-us\lambda}\tilde{\mu}(du)$. Since
$ng''(1-s\lambda/\sqrt{n})$ is bounded, (D2) holds and
$\lim_{n\rightarrow\infty}G_n(\lambda)=-(\tilde{\mu}(\{0\})+a)\lambda^2/2-\int_{(0,\infty)}
(e^{-\lambda u}-1+\lambda u)/u^2\tilde{\mu}(du)$. It follows in a
similar way that (E2) also holds and $\lim_{n\rightarrow\infty}
H_n(\lambda)=-(\tilde{\nu}(\{0\})+\omega^2-\omega)\lambda^2/2$
$-\int_{(0,\infty)}(e^{-\lambda u}-1+\lambda u)/u^2\tilde{\nu}(du)$.
\qed
 \bexample\label{example of diffusion}(\cite[Theorem
2.2]{I05})\;Assume that (D1), (E1) and the following conditions
hold:
\begin{itemize}
\item[(b.1)]\;
 $n\pi_n\rightarrow \pi$ and  $r_n\rightarrow r$ as $n\rightarrow\infty$
 for some
 $\pi\geq0$ and $r\geq0$,
\item[(b.2)]\;
 $n\bfE\Big[(\xi_n(1,1)-m_n)^2{\bf1}_{\{|\xi_n(1,1)-m_n|>\sqrt{n}\,\varepsilon}\Big]
 \rightarrow0$ as $ n\rightarrow\infty$ for all
 $\varepsilon>0$,
\item[(b.3)]\;
$\bfE\Big[(\eta_n(1)-\omega_n)^2{\bf
1}_{\{|\eta_n(1)-\omega_n|>\sqrt{n}\,\varepsilon\}}\Big]\rightarrow0$
as $n\rightarrow\infty$ for all $\varepsilon>0$.
  \end{itemize}
In this case, we see that conditions (a.1,2) are satisfied with
$\tilde{\mu}(du)=\pi\delta_0(du)$ and
$\tilde{\nu}(du)=r\delta_0(du)$ ($\delta_{x}(du)$ denote the dirac
measure at $u=x$). Then we still have the above limit theorem, and
the fluctuation limit process $Z(\cdot)$ is given by
 \beqlb
 dZ(t)=aZ(t)dt + \sqrt{\varrho(t)}\,dB(t),\qquad Z(0)=0,
 \label{OU diffusion}
 \eeqlb
where $\varrho(t)=r+\pi\phi(t)$, and $B(t)$ is a one-dimensional
Brownian motion.  \eexample

  \bexample\label{jump processes}\; Suppose
that $\bfP(\xi_n(1,1)=[\sqrt{n}\,])=1/n^2$ and
$\bfP(\xi_n(1,1)=1)=1-1/n^2$, while
$\bfP(\eta_n(1)=[\sqrt{n}\,])=1/n$ and $\bfP(\eta_n(1)=1)=1-1/n$. We
see that $n(m_n-1)\rightarrow0$, $\omega_n\rightarrow1$ and (a.1,2)
are satisfied with $\tilde{\mu}(du)=\tilde{\nu}(du)=\delta_1(du)$.
Another example is as follows. Suppose that $\tilde{\mu}(du)$ is any
non-degenerate finite measure on $(0,\infty)$. For large enough $n$,
let $\mu(du)=u^{-2}\tilde{\mu}(du)$,
$\mu_n(du)=\mu((1/n^{1/4},\infty))^{-1}1_{\{u>n^{1/4}\}}\mu(du)$,
and
  \beqnn
  g_n(s)=p_n\int_0^\infty e^{-\sqrt{n}u(1-s)}\mu_n(du)+
  (1-p_n),\quad\mbox{where}\quad
  p_n=\frac{\mu((1/n^{1/4},\infty))}{\mu((1/n^{1/4},\infty))+n^2}.
  \eeqnn
 Let $\xi_n(1,1)$ have the distribution corresponding to
$g_n(\cdot)$. Note that $\int_{\{u>1/n^{1/4}\}}u\mu(du)/\sqrt{n}
\rightarrow0$ and $\mu((1/n^{1/4},\infty))/n\rightarrow0$. Then it
is not hard to see that $n(m_n-1)\rightarrow0$ and condition (a.1)
is fulfilled with $\tilde{\mu}(du)$. $\eta_n(1)$ can be constructed
in a similar way. \eexample

 \bcorollary\label{example2} Let $L(x)$ and $L^*(x)$ be
positive functions slowly varying at $\infty$ such that $L(x)\sim
L^*(x)$ as $x\rightarrow\infty$. Consider a sequence of
GWI-processes with
 $(g_n, h_n)$ given by
  \beqnn
  g_n(s) \ar=\ar
  (1-m_n)+m_ns+\frac{\gamma}{n}(1-s)^{\alpha}L\Big(\frac{1}{1-s}\Big),\\
  h_n(s) \ar=\ar
  (1-\omega_n)+\omega_ns+\varpi(1-s)^{\alpha}L^*\Big(\frac{1}{1-s}\Big),
  \eeqnn
 where $1<\alpha\leq2$, $\gamma>0$, $\varpi>0$. $m_n$ and $\omega_n$
satisfy conditions (D1) and (E1).  Let $Z_n(\cdot)$ be defined by
(\ref{Z_n}) with $c_n$ satisfying
 \beqlb
 c_n\sim [nL(c_n)]^{1/\alpha}\quad (\sim[nL^*(c_n)]^{1/\alpha}).
 \label{c_nL}
 \eeqlb
Then $Z_n(\cdot)$ converges in distribution on $D([0,\infty),
\mathbb{R})$ to a OU type process $Z(\cdot)$ defined by
 \beqlb
 dZ(t)=aZ(t)dt + \sqrt[\alpha]{\varrho_1(t)}\,dX(t), \qquad
 Z(0)=0,
 \label{stable}
 \eeqlb
where $\varrho_1(t)=\varpi+\gamma\phi(t)$,
$\phi(t)=\omega\int_0^te^{au}du$, and $X(t)$ is a spectrally
positive $\alpha$-stable L\'{e}vy process with Laplace exponent
$-\lambda^{\alpha}$.
 \ecorollary

\proof By (\ref{c_nL}), $c_n/n\sim c_n^{1-\alpha}L(c_n)\rightarrow0$
and $c_n^2/n\sim c_n^{2-\alpha}L(c_n)\rightarrow\infty$, as
$n\rightarrow\infty$. Then (C) holds. Without loss of generality, we
consider $G_n(\lambda)$ on $\lambda\in[0,1]$ (see (\ref{G_n})). For
$\lambda>0$,
$G_n(\lambda)=-\gamma\lambda^\alpha[nL(c_n)/c_n^\alpha][L(c_n/\lambda)/L(c_n)]$,
and thus
 $
\lim_{n\rightarrow\infty}G_n(\lambda)=-\gamma\lambda^{\alpha}.
 $
For $\lambda=0$, the limit is trivial. Furthermore we have
$\sup_n|G'_n(\lambda)|$ is bounded in $\lambda\in(0,1]$ as in the
proof of Corollary \ref{example(1)}. Note that $G'_n(0)=0$ and thus
(D2) holds. It follows in a similar way that (E2) also holds and
 $
 \lim_{n\rightarrow\infty}H_n(\lambda)=-\varpi\lambda^{\alpha}.
 $
By Theorem \ref{OU limit}, the limit process $Z(\cdot)$ is described
by (\ref{OU equation}) with
$\mu(du)=(\gamma\alpha(\alpha-1)/\Gamma(2-\alpha))u^{-1-\alpha}du$,
$\nu(du)=(\varpi\alpha(\alpha-1)/\Gamma(2-\alpha))u^{-1-\alpha}du$,
$\varrho(\cdot)\equiv0$ and $\beta_1=\beta_2=0$. Note that
$\varrho_1(t)>0$ for $t\geq0$. Define the process
 \beqlb
X(t)=\int_0^t\int_{\mathbb{R}_+}\varrho_1(s)^{-\frac{1}{\alpha}}u\tilde{N}_0(ds,du)
+\int_0^t\int_{\mathbb{R}_+}\int_0^{\infty}\varrho_1(s)^{-\frac{1}{\alpha}}u{\bf
1}_{(0,\phi(s)]}(\zeta)\tilde{N}_1(ds,du,d\zeta).
 \label{stable construction}
 \eeqlb
Then $X(\cdot)$ is a martingale. By It\^{o}'s formula, it is not
hard to show that $X(\cdot)$ is a one-sided $\alpha$-stable process
with Laplace exponent $-\lambda^{\alpha}$. Thus we have
(\ref{stable}) by (\ref{stable construction}) and (\ref{OU
equation}).\qed

\section{Asymptotic results for estimators}

\setcounter{equation}{0}

In this section, we consider the statistical applications of our
limit theorems as in \cite{S94,I05}. For $n\in\mathbb{N}$, suppose
that a sequence of samples $\{(y_n(k),\eta_n(k)), k=1,2,\cdots, n\}$
is available. Then a natural estimator of the offspring mean $m_n$
is given by
 \beqnn
 \check{m}_n=\Big[\sum_{k=1}^ny_n(k-1)\Big]^{-1}\sum_{k=1}^n(y_n(k)-\eta_n(k)).
 \eeqnn
Using Theorem \ref{Theorem 1}, we can derive the following
asymptotic result for $\check{m}_n$.

\btheorem\label{asymptotic2.3}\;If the conditions of Theorem
\ref{Theorem 1} are fulfilled with $F(\lambda)$ being a unbounded
function, and $R'_n(0)\rightarrow a$ as $n\rightarrow\infty$ for
some $a\in\mathbb{R}$, then
  \beqlb
  n(\check{m}_n-m_n)\overset{d}{\longrightarrow}
  \frac{Y(1)-a\int_0^1Y(s)ds-Y'(1)}{\int_0^1Y(t)dt},
   \label{nature}
    \eeqlb
where $Y(t)$ and $Y'(t)$ are defined in Theorem \ref{Theorem 1}.
 \etheorem
Obviously, the above theorem applies to the case of Corollary
\ref{example(1)}. Compared with the result of \cite[Corollary
3.1]{S94}, this implies that the heavy-tailed stable distributions
of offspring and immigration variables do not affect the rate of
convergence of $\check{m}_n$ in the critical GWI-process.

We are also interested in the case when the conditions of Theorem
\ref{OU limit} are satisfied. Then $Y_n(\cdot)/n$ converges weakly
to the deterministic function $\phi(t)=\omega\int_0^te^{as}ds$,
which implies only that
$n(\check{m}_n-m_n)\overset{p}{\longrightarrow}0$ by Theorem
\ref{asymptotic2.3}. Thus we further consider the applications of
our fluctuation limit theorem, related to the CLSE of the offspring
mean $m_n$ based on only the information on $\{y_n(k)\}$ as follows.
For $n,k\geq1$ let $\mathcal{F}^n_k$ denote the $\sigma$-algebra
generated by $\{y_n(j), j=0,1,\cdots,k\}$. From
(\ref{GWI-equationn}),
 \beqlb
 \bfE[y_n(k)|\mathcal{F}^n_{k-1}]=m_ny_n(k-1)+\omega_n,\qquad
 n\geq1.
 \label{regression equation}
 \eeqlb
If we assume that the immigration mean $\omega_n$ is known, then the
CLSE $\hat{m}_n$ of $m_n$, based on (\ref{regression equation}), is
given by
 \beqlb
 \hat{m}_n=\frac{\sum_{k=1}^ny_n(k-1)(y_n(k)-\omega_n)}
 {\sum_{k=1}^ny_n^2(k-1)}.
 \label{form}
 \eeqlb
If $\omega_n$ is unknown, it is not hard to see that the joint CLSE
$(\tilde{m}_n,\tilde{\omega}_n)$ of $(m_n, \omega_n)$ is given by
 \beqnn
 \tilde{m}_n=\frac{\sum_{k=1}^ny_n(k-1)(y_n(k)-\overline{y_n}\,)}
 {\sum_{k=1}^n(y_n(k-1)-\overline{y^*_n}\,)^2},\qquad
 \tilde{\omega}_n=\overline{y_n}-\tilde{m}_n\overline{y^*_n},
  \eeqnn
where
  \beqlb
  \overline{y_n}=\frac{1}{n}\sum_{k=1}^ny_n(k),\qquad
  \overline{y^*_n}=\frac{1}{n}\sum_{k=1}^ny_n(k-1).
  \label{y_n*}
  \eeqlb
Using Theorem \ref{OU limit}, we can derive the following asymptotic
result for $\hat{m}_n$, $\tilde{m}_n$, and $\tilde{\omega}_n$, which
generalizes the result of \cite[Theorem 3.1]{I05}.
 \btheorem\label{asymptotic}\;If the conditions of Theorem \ref{OU
 limit} are fulfilled with $\omega>0$, then
  \beqlb
  \frac{n^2}{c_n}\,(\hat{m}_n-m_n)\overset{d}{\longrightarrow}
  \frac{\int_0^1\phi(t)d{M}(t)}{\int_0^1\phi^2(t)dt},
  \label{asymptotic1}
  \eeqlb
where $\phi(t)=\omega\int_0^te^{as}ds$, $Z(t)$ is defined by
(\ref{OU}), ${M}(t)= Z(t)-\int_0^taZ(s)ds$ and it can be regarded as
a deterministically-time-changed L\'{e}vy process. Furthermore,
 \begin{eqnarray}
\left(\begin{array}{c}
\frac{n^2}{c_n}(\tilde{m}_n-m_n)  \\
\\
\frac{n}{c_n}(\tilde{\omega}_n-\omega_n)
\end{array}\right)\overset{d}{\longrightarrow}
\left(\begin{array}{c}
\frac{\int^1_0\phi(t)d{M}(t)-{M}(1)\int_0^1\phi(t)dt}
{\int_0^1\phi^2(t)dt-(\int_0^1\phi(t)dt)^2}  \\
 \\
\frac{{M}(1)\int_0^1\phi^2(t)dt-\int_0^1\phi(t)dt\int^1_0\phi(t)d{M}(t)}
{\int_0^1\phi^2(t)dt-(\int_0^1\phi(t)dt)^2}
\end{array}\right).
 \label{asymptotic1+}
 \end{eqnarray}
   \etheorem

\bcorollary\label{example3}\;Consider the CLSE $\hat{m}_n$ of $m_n$
in the case of Corollary \ref{example2} when condition (E1) holds
with $\omega>0$. By the above theorem, we have
 \beqnn
 \frac{n^2}{c_n}\,(\hat{m}_n-m_n)\overset{d}{\longrightarrow}
 U:=
\frac{\int_0^1\phi(t)\sqrt[\alpha]{\varrho_1(t)}dX(t)}
  {\int_0^1\phi^2(t)dt},
   \eeqnn
where $1<\alpha\leq2$, $c_n$ is given by (\ref{c_nL}), $\phi(t)$,
$\varrho_1(t)$ and $X(t)$ are given in (\ref{stable}). It is easy to
see that $U$ has a $\alpha$-stable distribution and its Laplace
transform equals
 \beqnn
 \bfE[e^{-\lambda
 U}]=\exp\bigg\{\frac{\int_0^1\phi^{\alpha}(t)\varrho_1(t)dt}{(\int_0^1
 \phi^2(t)dt)^\alpha}\,\lambda^\alpha\bigg\},\qquad \lambda\geq0.
 \eeqnn
 \ecorollary
Finally we turn to the case when the conditions of Corollary
\ref{example1} are satisfied. In this case, it is possible to
consider the CLSE estimates for the offspring and immigration
variances $\pi_n$ and $r_n$. Let
$u_n(k)=y_n(k)-m_ny_n(k-1)-\omega_n$. Note that
 \beqlb
 \bfE[u^2_n(k)|\mathcal{F}^n_{k-1}]=\pi_ny_n(k-1)+r_n,\qquad
 n\geq1.
 \label{variance equation}
 \eeqlb
As in \cite{W91}, if we suppose that $m_n$ and $\omega_n$ are known,
then the joint CLSE $(\hat{\pi}_n, \hat{r}_n)$ of $(\pi_n, r_n)$,
based on (\ref{variance equation}), is given by
 \beqlb
 \hat{\pi}_n=\frac{\sum_{k=1}^nu_n^2(k)(y_n(k-1)-\overline{y^*_n}\,)}
 {\sum_{k=1}^n(y_n(k-1)-\overline{y^*_n}\,)^2},\qquad
 \hat{r}_n=\sum_{k=1}^nu^2_n(k)/n-\hat{\pi}_n\overline{y^*_n},
  \label{estimator variance}
  \eeqlb
where $\,\overline{y^*_n}$ is defined by (\ref{form}). If $m_n$ and
$\omega_n$ are unknown, we can use
$\hat{u}_n(k)=y_n(k)-\hat{m}_ny_n(k-1)-\hat{\omega}_n$ instead of
$u_n(k)$ in (\ref{estimator variance}) and we get another joint CLSE
denoted by $(\tilde{\pi}_n,\tilde{r}_n)$. Using Theorem \ref{OU
 limit} again, we have the following asymptotic result for the above estimators, where the
 jumps of the fluctuation limit obviously play an important role.

\btheorem\label{asymptotic variance}\;If the conditions of Corollary
\ref{example1}, i.e. (D1), (E1) and (a.1,2) are fulfilled  with
$\omega>0$, then
 \begin{eqnarray}
\left(\begin{array}{c}
n(\hat{\pi}_n-\pi_n)  \\
\\
\hat{r}_n-r_n
\end{array}\right)\overset{d}{\longrightarrow}
\left(\begin{array}{c}
\frac{\int^1_0\phi(t)d{J}(t)-{J}(1)\int_0^1\phi(t)dt}
{\int_0^1\phi^2(t)dt-(\int_0^1\phi(t)dt)^2}  \\
 \\
\frac{{J}(1)\int_0^1\phi^2(t)dt-\int_0^1\phi(t)dt\int^1_0\phi(t)d{J}(t)}
{\int_0^1\phi^2(t)dt-(\int_0^1\phi(t)dt)^2}
\end{array}\right),
 \label{asymptotic2+}
 \end{eqnarray}
 where $\phi(t)=\omega\int_0^te^{as}ds$, and $J(t)$ is a martingale
defined by
  \beqlb
  J(t)=\int_0^t\int_{\mathbb{R}_+}
  u^2\tilde{N}_0(ds,du)+\int_0^t\int_{\mathbb{R}_+}\int_0^{\phi(s)}
  u^2\tilde{N}_1(ds,du,d\zeta).
   \label{J(t)}
  \eeqlb
Here $\tilde{N}_0$ and $\tilde{N}_1$ are the compensated Poisson
random measures given in (\ref{OU_square equation}). Moreover,
(\ref{asymptotic2+}) still holds if $\hat{\pi}_n$ and $\hat{r}_n$
are replaced by $\tilde{\pi}_n$ and $\tilde{r}_n$.
 \etheorem

\bremark\label{remark 3.1}\; $N_0$ and $N_1$ are the Poisson random
measures given in (\ref{OU_square equation}) with intensities
$ds\nu(du)$ and $ds\mu(du)d\zeta$, and $\int_0^\infty
u^2\nu(du)+\int_0^\infty u^2\mu(du)<\infty$. Let the limiting random
vector in (\ref{asymptotic2+}) be denoted by $(U_1,U_2)^T$. It is
not hard to see that if $\int_0^\infty u^4\nu(du)+\int_0^\infty
u^4\mu(du)<\infty$, then
 \beqnn
\bfE[U_1^2]
 &=&
L^{-2}\bigg[\int_0^1\bigg(\phi(t) - \int_0^1\phi(s)d{s\bigg)^2} dt
\int_0^\infty u^4\nu(du) \\
 & &\qquad\qquad
+ \int_0^1\phi(t)\bigg(\phi(t) - \int_0^1\phi(s)d{s\bigg)^2} dt
\int_0^\infty u^4\mu(du)\bigg],
 \eeqnn
where $L=\int_0^1\phi^2(t)dt-(\int_0^1\phi(t)dt)^2$. Otherwise we
have $\bfE[U_1^2]=\infty$. Note that $\phi(\cdot)$ is not a const
function. So if $\nu\neq0$ or $\mu\neq0$ (equivalently $N_0$ or
$N_1$ is not degenerate), then $U_1$, and similarly $U_2$, are not
degenerate.

We see that if the conditions of Corollary \ref{example1} are
fulfilled with $\tilde{\nu}(\mathbb{R}_+\backslash\{0\})>0$ or
$\tilde{\mu}(\mathbb{R}_+\backslash\{0\})>0$, which means that the
sequence of the offspring (or immigration) distributions fails to
satisfy Lindeberg condition (b.2) (or (b.3)), then the resulting
fluctuation limit $Z(\cdot)$ is a OU type process with positive
jumps (see (\ref{OU_square equation})). Thus, in this case,
$\hat{\pi}_n$ has the limit law $U_1$ with normalizing factor $n$,
and $\hat{r}_n$ is not a consistent estimator. However, if we return
to the case of \cite{I05} (see Example \ref{example of diffusion}),
which implies that Lindeberg conditions are satisfied and the
resulting fluctuation limit $Z(\cdot)$ is a OU diffusion process
without jumps (see (\ref{OU diffusion})), then
$n(\hat{\pi}_n-\pi_n)\overset{p}{\rightarrow}0$ and
$\hat{r}_n-r_n\overset{p}{\rightarrow}0$. In this case, to get the
appropriate rates of convergence for $\hat{\pi}_n$ and $\hat{r}_n$,
we give the following theorem.
 \eremark
\btheorem\label{diffusion variance theorem}\; Consider the case of
Example \ref{example of diffusion}. Let $a_{4,n} = \bfE[(\xi_n(1,1)
- m_n)^4]$ and $b_{4,n} = \bfE[(\eta_n(1) - \omega_n)^4]$. Suppose
that (D1), (E1), (b.1) and the following conditions hold with
$\omega>0$:
 \begin{itemize}

 \item[{(c.1)}]\;
$na_{4,n}\rightarrow a_4$ and $b_{4,n}\rightarrow b_4$ as
$n\rightarrow \infty$ for some $a_4\geq0$ and $b_4\geq0$,

 \item[{(c.2)}]\;
$n\bfE\Big[(\xi_n(1,1)-m_n)^4{\bf1}_{\{(\xi_n(1,1)-m_n)^2>\sqrt{n}\,
\varepsilon \}}\Big]\rightarrow0$ as $ n\rightarrow\infty$ for all
$\varepsilon>0$,

\item[(c.3)]\;
$\bfE\Big[(\eta_n(1)-\omega_n)^4{\bf 1}_{\{(\eta_n(1)-\omega_n)^2>
\sqrt{n}\, \varepsilon\}}\Big]\rightarrow0$ as $n\rightarrow\infty$
for all $\varepsilon>0$.

  \end{itemize}
Let $\phi(t)=\omega\int_0^te^{as}ds$, $\varrho_2(t) =
2\pi^2\phi^2(t) + (a_4+4\pi r)\phi(t)+(b_4-r^2)$ and $V(t) =
\int_0^t \sqrt{\varrho_2(s)}\, dW(s)$, where $W(t)$ is a
one-dimensional Brownian motion. Let $\Sigma =
\big(\int_0^1\phi^2(t)dt - (\int_0^1
\phi(t)dt)^2\big)^{-2}(\sigma_{ij})_{2\times 2}$, where
   \beqnn
\left.\begin{array}{ll} \quad\sigma_{11}=\int_0^1\big(\phi(t) -
\int_0^1\phi(s)ds\big)^2 \varrho_2(t)dt,\quad
\sigma_{22}=\int_0^1\big(\int_0^1\phi^2(s)ds-\phi(t)\int_0^1
\phi(s)ds\big)^2\varrho_2(t)dt,\\
\\
\quad\sigma_{12}=\sigma_{21}=\int_0^1\big(\phi(t)-\int_0^1\phi(s)ds\big)
\big(\int_0^1\phi^2(s)ds-\phi(t)\int_0^1\phi(s)ds\big) \varrho_2(t)dt. \\
  \end{array}
  \right.
  \eeqnn
Then we have
 \begin{eqnarray}
\left(\begin{array}{c}
n^{3/2}(\hat{\pi}_n-\pi_n)  \\
\\
n^{1/2}(\hat{r}_n-r_n)
\end{array}\right)\overset{d}{\longrightarrow}
\left(\begin{array}{c}
\frac{\int^1_0\phi(t)d{V}(t)-{V}(1)\int_0^1\phi(t)dt}
{\int_0^1\phi^2(t)dt-(\int_0^1\phi(t)dt)^2}  \\
 \\
\frac{{V}(1)\int_0^1\phi^2(t)dt-\int_0^1\phi(t)dt\int^1_0\phi(t)d{V}(t)}
{\int_0^1\phi^2(t)dt-(\int_0^1\phi(t)dt)^2}
\end{array}\right)\overset{d}{=}\mathcal{N}(0,\Sigma).
 \label{asymptotic2+2}
 \end{eqnarray}
Furthermore, (\ref{asymptotic2+2}) still holds if $\hat{\pi}_n$ and
$\hat{r}_n$ are replaced by $\tilde{\pi}_n$ and $\tilde{r}_n$,
respectively.
  \etheorem
\bremark It is easy to see that the above condition (c.1) implies
that (b.2) and (b.3) in Example \ref{example of diffusion} hold. So
the conditions of our theorem is in the case of Example \ref{example
of diffusion}. If either of $\pi$, $a_4$ and $b_4-r^2$ is not $0$,
then the limit normal law in (\ref{asymptotic2+2}) is not
degenerate. \eremark

\bexample\label{example of diffusion variance}(\cite[Example
2.1]{I05}) The conditions of Theorem \ref{diffusion variance
theorem} are satisfied for the following examples with
$n(m_n-1)\rightarrow a$, $n\pi_n\rightarrow a$, $na_{4,n}\rightarrow
a$ and $n\bfE[(\xi_n(1,1)-m_n)^6]\rightarrow a$, as
$n\rightarrow\infty$ for some $a\geq0$. {\rm (i)}\, $\xi_n(1,1)$ has
a Bernoulli distribution with mean $1-an^{-1}$. {\rm (ii)}\, the
offspring distributions are geometric distributions with parameter
$p_n=1-an^{-1}$, i.e. $\bfP(\xi_n(1,1)=i)=p_n(1-p_n)^{i-1}$,
$i=1,2,\cdots$. \eexample

\section{Proof of main results}

\setcounter{equation}{0}

{\it Proof of Theorem \ref{Theorem 1}\;} For the proofs of Lemma
\ref{prop2.1} and \ref {prop2.2}, we can follow the proof of Lemma
\ref{lemma2.3} or apply directly \cite[Corollary 1,2]{L91}. So we
skip them. Now the limit functions $R$ and $F$ have representations
(\ref{R}) and (\ref{F}).  Fix $0\leq\lambda\leq M$ for any constant
$M>0$. Let $\lambda_n=b_n(1-e^{-\lambda/b_n})$ and we have
$\lambda_n\rightarrow\lambda$. It follows from condition (A) that
$|R_n(\lambda_n)-R_n(\lambda)|\leq k(M)|\lambda_n-\lambda|$, where
$k(M)>0$ is a constant, and that
$\lim_{n\rightarrow\infty}R_n(\lambda_n)=R(\lambda)$. By condition
(B) and the fact that $F_n$ is a nondecreasing function on
$\lambda\in[0,M]$ for sufficiently large $n$, we have
$F_n\rightarrow F$ locally uniformly. It implies that
$\lim_{n\rightarrow\infty}F_n(\lambda_n)=F(\lambda)$. Let
$\tilde{R}_n(\lambda)=R_n(\lambda_n)$ and
$\tilde{F}_n(\lambda)=F_n(\lambda)$. Note that the sequence
$\{(Y_n(\frac{l}{n}), Y_n'(\frac{l}{n})), l\in\mathbb{N}\}$ is a
Markov chain with state space $\hat{E}_n:=\{(i/b_n,j/b_n):
(i,j)\in\mathbb{N}^2\}$ and the (discrete) generator $A_n$ of
$\{(Y_n(t), Y_n'(t)), t\geq0\}$ is given by \beqnn
 A_n
 e^{-\langle z,x\rangle}
 &=&
 n\Big[\big(g_n(e^{-z_1/b_n})\big)^{b_nx_1}
  h_n(e^{-(z_1+z_2)/b_n})e^{-z_2x_2}-e^{-\langle z,x\rangle}\Big]\\
 &=&
e^{-\langle
z,x\rangle}n\Big[\exp\big\{-x_1\alpha_n(z)e^{z_1/b_n}\tilde{R}_n(z_1)/n\big\}
\exp\big\{-\beta_n(z)\tilde{F}_n(z_1+z_2)/n\big\}-1\Big]\\
 &=&
 -e^{-\langle z, x\rangle}
 \big[x_1\alpha_n(z)e^{z_1/b_n}\tilde{R}_n(z_1)+
 \beta_n(z)\tilde{F}_n(z_1+z_2)\big] + o(1),
 \eeqnn
where $x\in\hat{E}_n$, $z=(z_1,z_2)\gg0$, $\alpha_n(z)
=\big(e^{z_1/b_n}g_n(e^{-z_1/b_n})-1\big)^{-1}
\ln\big(e^{z_1/b_n}g_n(e^{-z_1/b_n})\big)$, and
$\beta_n(z)=\big(h_n(e^{-(z_1+z_2)/b_n})-1\big)^{-1}
 \ln\big(h_n(e^{-(z_1+z_2)/b_n})\big)$. On the other hand, let $A$ be the
infinitesimal generator of $(Y(\cdot),Y'(\cdot))$. For $z\gg0$ and
$x\in\mathbb{R}_+^2$,
 \beqnn
 Ae^{-\langle z,x\rangle}=-e^{-\langle z, x\rangle}\big[x_1R(z_1)+F(z_1+z_2)\big].
 \eeqnn
We need to prove that
 $
 \lim_{n\rightarrow\infty}\sup_{x\in
 \hat{E}_n}|A_ne^{-\langle z,x\rangle}-Ae^{-\langle z,x\rangle}|=0.
 $
The remaining proof is essentially the same as that of in
\cite[Theorem 2.1]{L05} or \cite[Theorem 2.1]{M09} and so we omit
it. \qed

 Let us write $f\in C_{\ast}(\mathbb{R})$ if $f$ is a bounded
continuous function from $\mathbb{R}$ to $\mathbb{R}$ satisfying
$f(x) = o(x^2)$ when $x\rightarrow0$. Let $\Gamma=[-1,\infty)$ and
$\Gamma_n=\{ (i-1)/c_n: i\in\mathbb{N}\}$. Let $\mu_n$ be the
distribution of $\frac{\xi_n(1,1)-1}{c_n}$. Then for sufficiently
large $n$, $\mu_n$ is a probability measure on $\Gamma$ supported by
$\Gamma_n$.
\\
\\
{\it Proof of Lemma \ref{lemma2.3}} (sketch)\; Set $S_n(\lambda)
 =
n^2\big[e^{-\lambda/c_n}\big(1-(m_n-1)\lambda/c_n\big)-g_n\big(e^{-\lambda/c_n}\big)\big]$
and it follows from mean-value theorem that
 \beqlb
S_n(\lambda)
 \ar=\ar
G_n(\lambda) + n^2[m_n-g'_n(\vartheta_n)](e^{-\lambda/c_n}
-1+\lambda/c_n) \nonumber\\
 \ar \ar
+n^2(1-m_n)(e^{-\lambda/c_n}
-1+\lambda/c_n)+n^2(m_n-1)(1-e^{-\lambda/c_n})\lambda/c_n,
 \label{equ2.1}
 \eeqlb
where $1-\lambda/c_n\leq \vartheta_n\leq e^{-\lambda/c_n}$. Under
condition (D2), the sequence
 $
|G'_{n}(\lambda)|=n^2|g'_{n}(1-\lambda/c_n)-m_n|\big/c_n
 $
is uniformly bounded on each bounded interval $[0,c]$ for $c\geq0$
and thus the sequence $n^2|g'_n(\vartheta_n)-m_n|\big/c_n$ is also
uniformly bounded. By (C), (D1) and (D2), we have
 $
S_n(\lambda)\rightarrow G(\lambda)+\frac{1}{2}a\gamma_0,
 $
as $n\rightarrow\infty$. To get (\ref{G}), it is enough to consider
the limit representation of $S_n$. Note that
 \beqnn
 e^{\lambda/c_n}S_n(\lambda)=-n^2\int_{\Gamma}(e^{-\lambda u}-1+\lambda
 u)\mu_n(du).
 \eeqnn
We can use Venttsel's classical method (see \cite{V59}) to prove it.
More precisely, by modifying slightly the proofs of Proposition 2.1
and 3.1 in \cite{LM08}, we can show that there exist some constants
$\hat{\beta}_1\in\mathbb{R}$, $\hat{\sigma}_1\geq0$, and a
$\sigma$-finite measure $\mu$ defined as in (\ref{G}) such that
 \begin{itemize}
\item[\rm{(i)}]\;
 $\displaystyle n^2\int_{\Gamma} (\chi(u)-u)
\mu_n(du) \rightarrow\hat{\beta}_{1}
 $\quad
as $n\rightarrow \infty$;

\item[\rm{(ii)}]\;
 $\displaystyle n^2\int_{\Gamma} \chi^2(u) \mu_n(d
u) \rightarrow2\hat{\sigma}_1 + \int_0^\infty \chi^2(u)\mu(du)
 $\quad
as $ n\rightarrow\infty$;

\item[\rm{(iii)}]\;
 $
\displaystyle\lim\limits_{n\rightarrow\infty} n^2\int_{\Gamma} f(u)
\mu_n(du) = \int_0^\infty f(u) \mu(du),
 $\quad for $f\in C_{\ast}(\mathbb{R})$.
  \end{itemize}
Note that $e^{-\lambda x} -1+ \lambda
\chi(x)-\frac{1}{2}\lambda^2\chi^2(x)\in C_{\ast}(\Gamma)$ as a
function of $x\in\Gamma$ for fixed $\lambda\geq0$. The above results
imply that the limit function of $S_n$ has a L\'{e}vy-Khintchine
type representation. Let
$\sigma_1=\hat{\sigma}_1+\frac{1}{2}a\gamma_0$ and let
$\beta_1=\hat{\beta}_1+\int_0^\infty(u-\chi(u))\mu(du)$. Then we
have (\ref{G}). But we still need to verify $\sigma_1\geq0$. It
follows from {\rm (a.1)}, (C) and (D1) that
 \beqlb
\frac{n^2}{c_n}\int_{\Gamma} \chi(u) \mu_n(du)
=\frac{n^2}{c_n}\int_{\Gamma} (\chi(u)-u) \mu_n(du) +
\frac{n^2(m_n-1)}{c^2_n},
 \label{cutmoment}
 \eeqlb
which tends to $a\gamma_0$ as $n\rightarrow\infty$. Let $E$ be the
set of $\varepsilon>0$ for which $\mu( |u|=\varepsilon)=0$. By
(\ref{cutmoment}), {\rm (ii)} and {\rm (iii)}, we obtain
 \beqnn
\lim_{E\ni\varepsilon\downarrow0}
\lim_{n\rightarrow\infty}n^2\int_{\{|u|<\varepsilon\}}
\Big(\chi^2(u)+\frac{\chi(u)}{c_n}\Big)\mu_n(du)=2\sigma_1.
 \eeqnn
The support of $\mu_n$ is $\Gamma_n$ and for large enough $n$,
$\chi^2(u)+\big(\chi(u)\big/c_n\big)\geq0$ if $u\in\Gamma_n$. Thus
$\sigma_1\geq0$.\qed

Let $\hat{\Gamma}_n=\{ i/c_n: i\in\mathbb{N}\}$ and let $\nu_n$ be
the distribution of $\frac{\eta_n(1)}{c_n}$. Then $\nu_n$ is a
probability measure on $[0,\infty)$ supported by $\hat{\Gamma}_n$.

 \blemma\label{lemma2.5} Under conditions (C), (E1) and (E2), (\ref{H})
holds. As $n\rightarrow\infty$, we also have
\begin{itemize}

\item[\rm{(i)}]\;
 $\displaystyle
n\int_0^\infty (\chi(u)-u) \nu_n(du) \rightarrow
\beta_2-\int_0^\infty (u - \chi(u))\nu(du);
 $

\item[\rm{(ii)}]\;
 $\displaystyle
n\int_0^\infty\chi^2(u)\nu_n(du)\rightarrow2\sigma_2 +
\omega\gamma_0 + \int_0^\infty\chi^2(u)\nu(du);
 $

\item[\rm{(iii)}]\;
 $\displaystyle
\lim_{n\rightarrow\infty}n\int_0^\infty f(u) \nu_n(du)=\int_0^\infty
f(u) \nu(du)$,\quad
 for $f \in C_{\ast}(\mathbb{R}_+)$.
 \end{itemize}
 \elemma
\proof This lemma is proved with the same method as Lemma
\ref{lemma2.3}. But we need to prove that
$2\sigma_2+\omega\gamma_0\geq\omega^2\gamma_0$. Let
$\hat{a}_n=\int_0^\infty\chi(u)\nu_n(du)$ and let $\hat{E}$ be the
set of $\varepsilon>0$ for which $\nu(u =\varepsilon)=0$. By (C),
(E1), {\rm (i)} and {\rm (ii)}, it is not hard to show that
 \beqnn
 &&\lim_{\hat{E}\ni\varepsilon\downarrow0}
\lim_{n\rightarrow\infty}n\int_{\{u<\varepsilon\}}\big(\chi^2(u)-\frac{\chi(u)}{c_n}\big)\nu_n(du)
=2\sigma_2, \\
 &&\lim_{\hat{E}\ni\varepsilon\downarrow0}
\lim_{n\rightarrow\infty}n\int_{\{u<\varepsilon\}}(\chi(u)-\hat{a}_n)^2\nu_n(du)
=2\sigma_2+\omega\gamma_0-\omega^2\gamma_0.
 \eeqnn
For large enough $n$, $\chi^2(u)-\big(\chi(u)\big/c_n\big)\geq0$ if
$u\in\hat{\Gamma}_n$. Then we are finished.
  \blemma\label{lemma2.6}\; Under
the conditions of Theorem \ref{OU limit}, we have for $t\geq0$,
 \beqlb
 \limsup_{n\rightarrow\infty}\frac{1}{n}\bfE\Big[\sup_{0\leq s\leq
 t}Y_n(s)\Big]\leq |a| \Phi(t)+\omega t,
 \label{supinequality}
 \eeqlb
 \beqlb
 \limsup_{n\rightarrow\infty}\bfE\Big[\sup_{0\leq s\leq t}|Z_n(s)|\Big]
 \leq \hat{M}(t)\big[1+(|a|+1)t\,\exp\{(|a|+1)t\}\big],
 \label{Z_ninequality}
  \eeqlb
where $\Phi(t)=\int_0^t\int_0^se^{(|a|+1)u}duds$,
$\hat{M}(t)=2K(\Phi(t)+t) + 4\sqrt{K(\Phi(t)+t)}$, and $K$ is a
positive constant defined as in (\ref{L}).

 \elemma
\proof Note that (\ref{GWI-equationn}) can be rewritten into the
following form:
 \beqlb
 y_n(l)=\sum_{k=1}^l\sum_{j=1}^{y_n(k-1)} (\xi_n(k,j)-1) +
 \sum_{k=1}^l\eta_n(k).
  \label{difference equation}
 \eeqlb
Let $\hat{\xi}_n(k,j)=(\xi_n(k,j)-1)\big/c_n$,
$w_n(k)=\sum_{j=1}^{y_n(k-1)}\big(\chi(\hat{\xi}_n(k,j))-\bfE[\chi(\hat{\xi}_n(k,j))]\big)$,
and $W_n(l)=\sum_{k=1}^lw_n(k)$. let $\tilde{\mathcal{F}}^n_k$
denote the $\sigma$-algebra generated by $\{(w_n(j),y_n(j)),
j=0,1,\cdots,k\}$. Since
$\bfE[w_n(k)|\tilde{\mathcal{F}}_{k-1}^n]=0$, $W_n([nt])$ is a
square integrable martingale, and the quadratic variation is
$\sum_{k=1}^{[nt]}w^2_n(k)$. On the other hand, it follows from
conditions (D1) and (E1) that
 \beqlb
\frac{1}{n^2}\sum_{k=1}^{[nt]}\bfE[y_n(k-1)]=\omega_n\int_0^{[nt]/n}\int_0^{[ns]/n}m_n^{[nu]}duds,
 \label{expectation limit}
 \eeqlb
which tends to $\omega\int_0^t\int_0^s e^{au}duds$, as
$n\rightarrow\infty$. Then applying Doob's inequality to martingale
terms in (\ref{difference equation}), we have for sufficiently large
$n$,
 \beqnn
\bfE\Big[\sup_{0\leq s\leq t}Y_n(s)\Big]
 \ar\leq\ar
c_n\bfE\Big[\sum_{k=1}^{[nt]}\sum_{j=1}^{y_n(k-1)}\big(\hat{\xi}_n(k,j)-\chi(\hat{\xi}_n(k,j))\big)\Big]
+2c_n\bfE^{\frac{1}{2}}\Big[W^2_n([nt])\Big]\\
 \ar \ar
+\,c_n\sum_{k=1}^{[nt]}\bfE\big[y_n(k-1)\big]\big|\bfE\big[\chi(\hat{\xi}_n(k,j))\big]\big| + n\omega_nt\\
  \ar\leq\ar
n^2c_n\Phi(t)\int_\Gamma(u-\chi(u))\mu_n(du) +
2c_n\Big(n^2\Phi(t)\,\bfvar\,\chi(\hat{\xi}_n(1,1))\Big)^{\frac{1}{2}}\\
 \ar \ar
+\, n^2c_n\Phi(t)\Big|\int_\Gamma\chi(u)\mu_n(du)\Big|+n\omega_nt.
 \eeqnn
By {\rm (i), (ii)}, (C), (D1) and (E1), we obtain
$nc_n\int_{\Gamma}\chi(u)\mu_n(du)\rightarrow a$ and then
(\ref{supinequality}) holds. By (\ref{difference equation}), the
 sequence $Z_n(\cdot)$ are given by
 \beqlb
 Z_n(t)=\sum_{k=1}^{[nt]}(m_n-1)Z_n\big(\frac{k-1}{n}\big)
 +\sum_{k=1}^{[nt]}\sum_{j=1}^{y_n(k-1)}\frac{\xi_n(k,j)-m_n}{c_n}
 +\sum_{k=1}^{[nt]}\frac{\eta_n(k)-\omega_n}{c_n}.
 \label{GW difference equation}
 \eeqlb
Let $\hat{\eta}_n(k)=\eta_n(k)/c_n$. By Doob's inequality, it is not
hard to see that for sufficiently large $n$,
 \beqlb
\bfE\Big[\sup_{0\leq s\leq t}|Z_n(s)|\Big]
 \ar\leq\ar n|m_n-1|\int_0^t\bfE[|Z_n(s)]ds +
 2n^2\Phi(t)\int_{\Gamma}(u-\chi(u))\mu_n(du)\nonumber\\
 \ar \ar
 \,+2\Big(n^2\Phi(t)\,\bfvar\,\chi(\hat{\xi}_n(1,1))\Big)^{\frac{1}{2}}
 + 2nt\int_0^\infty(u-\chi(u))\nu_n(du)\nonumber\\
 \ar \ar
 \,+ 2\Big(nt\,\bfvar\,\chi(\hat{\eta}_n(1))\Big)^{\frac{1}{2}}.
 \label{Z_ninequality2}
 \eeqlb
By Gronwall's inequality and standard stopping argument,  {\rm (i),
(ii)}, (C), (D1) and Lemma \ref{lemma2.5} implies
 \beqlb
 \bfE[|Z_n(t)|]\leq2\big\{K(\Phi(t)+t) +
 2\big(K(\Phi(t)+t)\big)^{\frac{1}{2}}\big
 \}\exp\{(|a|+1)t\},
 \label{L}
 \eeqlb
where $
K=\sup_n\big(n^2\int_{\Gamma}(u-\chi(u)+\chi^2(u))\mu_n(du)+n\int_0^\infty(u-\chi(u)+\chi^2(u))\nu_n(du)\big)$.
By the above inequality and (\ref{Z_ninequality2}), we obtain
(\ref{Z_ninequality}).\qed
 \blemma\label{lemma2.7} Let $\phi_n(t)=\bfE[Y_n(t)]/n$ for $t\geq0$. Under the
conditions of Theorem \ref{OU limit}, the sequence
$(Z_n(\cdot),\phi_n(\cdot))$ is tight in
$D([0,\infty),\mathbb{R}\times\mathbb{R}_+)$.
 \elemma
\proof By Lemma \ref{lemma2.6}, $C(t):= 1+
 \limsup_{n\rightarrow\infty}\big(\frac{1}{n}\bfE\big[\sup_{0\leq s\leq
 t}Y_n(s)\big]+\bfE\big[\sup_{0\leq s\leq
 t}|Z_n(s)|\big]\big)$ is a locally bounded
function of $t\geq0$. Then $Z_n(t)$ is a tight sequence of random
variables for every $t\geq0$. Now let $\{\tau_n\}$ be a sequence of
stopping times bounded by $T$ and let $\{\delta_n\}$ be a sequence
of positive constants such that $\delta_n\rightarrow0$ as
$n\rightarrow0$. By Doob's Optional Sampling Theorem, we obtain as
in the calculations in (\ref{Z_ninequality2}) that for sufficiently
large $n$,
 \beqnn
 \ar\ar\bfE\big[\big|Z_n(\tau_n+\delta_n)-Z_n(\tau_n)\big|\big]\nonumber\\
 \ar\leq\ar
 2K\int_0^{\frac{[n\delta_n]+1}{n}}\frac{1}{n}\bfE[y_n([n\tau_n]+[ns])]\,ds
 +(|a|+1)\int_0^{\frac{[n\delta_n]+1}{n}}\bfE[Z_n(\frac{[n\tau_n]+[ns]}{n})]\nonumber\\
 \ar \ar
 +\,\Big(K\int_0^{\frac{[n\delta_n]+1}{n}}\frac{1}{n}\bfE[y_n([n\tau_n]+[ns])]\,ds\Big)^{\frac{1}{2}}
 +2K\Big(\delta_n+\frac{1}{n}\Big)+\Big(K(\delta_n+\frac{1}{n})\Big)^{\frac{1}{2}}
 \nonumber\\
  \ar\leq\ar
 (2K+|a|+1)\int_0^{\delta_n+\frac{1}{n}}C(T+s)ds + \Big(K\int_0^{\delta_n+\frac{1}{n}}C(T+s)ds\Big)^{\frac{1}{2}}\nonumber\\
  \ar \ar
  +2K\Big(\delta_n+\frac{1}{n}\Big)+\Big(K(\delta_n+\frac{1}{n})\Big)^{\frac{1}{2}}.
  \label{Aldous}
 \eeqnn
Then $Z_n(\cdot)$ is tight in $D([0,\infty),\mathbb{R})$ by the
criterion of Aldous \cite{A78}. It is easy to see that $\phi_n(t)$
converges to $\phi(t):=\omega\int_0^te^{as}ds$ in distribution on
$D([0,\infty),\mathbb{R}_+)$. By Jocod and Schiryaev \cite[Corollary
3.33, P.317]{JS87}, $(Z_n(\cdot),\phi_n(\cdot))$ is tight in
$D([0,\infty),\mathbb{R}\times\mathbb{R}_+)$.\qed

Let $Z(\cdot)$ be any limit point of $Z_n(\cdot)$. Without loss of
generality, by Skorokhod's theorem, we can assume that on some
Skorokhod's space $(\Omega, \mathcal{F}, \mathcal{F}_t, {\mathbf
P})$,
$(Z_n(\cdot),\phi_n(\cdot))\overset{a.s.}{\longrightarrow}(Z(\cdot),\phi(\cdot))$
in the topology of $D([0,\infty),\mathbb{R}\times\mathbb{R}_+)$.
 \blemma\label{lemma2.8} For any fixed $\lambda\in\mathbb{R}$,
 \beqlb
 L(t)=e^{i\lambda Z(t)}-e^{i\lambda Z(0)}-\int_0^t e^{i\lambda
 Z(s)}A(Z(s), \phi(s), \lambda)ds
 \label{martingale problem}
 \eeqlb
is a complex-valued local $\mathcal{F}_t$-martingale. Here $i^2=-1$
and
 \beqnn
A(x_1,x_2,\lambda)=ia\lambda x_1 +
(a\gamma_0\lambda^2/2-G(-i\lambda))x_2
+\gamma_0(\omega^2-\omega)\lambda^2/2-H(-i\lambda),
 \eeqnn
where $G$ and $H$ are defined by (\ref{G}) and (\ref{H}),
respectively.
 \elemma
\proof  Define the stopping times
 \beqlb
\tau^b \ar = \ar \inf\{t\geq0: |Z(t)|\geq b\ \mbox{or}\ |Z(t-)|\geq
b\},
 \nonumber\\
\tau_n^b \ar = \ar \inf\{t\geq0: |Z_n(t)|\geq b\ \mbox{or}\
|Z_n(t-)|\geq b\}.
 \nonumber
 \eeqlb
Let $Z^b(t)=Z(t\wedge\tau^b)$, $Z_n^b(t)=Z_n(t\wedge\tau_n^b)$, and
analogously  $\phi^b(t)$, $\phi^b_n(t)$. It follows from  \cite
[Proposition 2.11, P.305]{JS87} that for all but countably many $b$,
 \beqnn
 \tau^b_n  \overset{a.s.}{\longrightarrow}  \tau^b \ \mbox{in}\
 \mathbb{R} \quad\mbox{and}\quad
 (Z_n^b(\cdot),\phi_n^b(\cdot)) \overset{a.s.}{\longrightarrow}
 (Z^b(\cdot),\phi^b(\cdot))
 \eeqnn
in the topology of $D([0,\infty),\mathbb{R}\times\mathbb{R}_+)$.
Define $\tau^b_n(t)=\tau^b_n\wedge t$ and $\tau^b(t)=\tau^b\wedge
t$. We claim that
 \beqlb
 \tau^b_n(\cdot)\overset{a.s.}{\longrightarrow}\tau^b(\cdot)\quad \mbox{in}\
 C([0,\infty), \mathbb{R}_+),\ \quad \mbox{as}\ n\rightarrow\infty.
 \label{stopping times}
 \eeqlb
In fact, since $0\leq
\tau^b_n(t+\varepsilon)-\tau^b_n(t)\leq\varepsilon$ for any
$t\geq0$, the criterion of Aldous yields tightness for
$\{\tau^a_n(\cdot),\ n\geq1\}$. On the other hand,
$\{Z_n(\frac{k}{n}):k\geq1\}$ is a time-inhomogeneous Markov chain.
For fixed $\lambda\in\mathbb{R}$,
 \beqnn
 L_n(l)=e^{i\lambda Z_n(\frac{l}{n})}-e^{i\lambda
 Z_n(0)}-\sum_{k=0}^{l-1}\Big(\bfE\big[e^{i\lambda
 Z_n(\frac{k+1}{n})}\big|\mathcal{F}^n_{k}\big]-e^{i\lambda
 Z_n(\frac{k}{n})}\Big)
 \eeqnn
is a complex-valued martingale. (\ref{trans matrix}) implies that
 \beqlb
 L_n([nt])= e^{i\lambda Z_n(t)}-e^{i\lambda
 Z_n(0)}-\int_0^{\frac{[nt]}{n}}e^{i\lambda
 Z_n(s)}n\big[A_n\big(Z_n(s),\phi_n(s),\lambda\big)-1\big]ds,
 \eeqlb
where $A_n(x_1,x_2,\lambda)=e^{-i\lambda/c_n\,
 (n(m_n-1)x_2+\omega_n)}\big(e^{-i\lambda/c_n}g_n(e^{i\lambda/c_n})\big)^{c_nx_1+nx_2}h_n(e^{i\lambda/c_n})$.
 For simplicity, we denote $L_n([nt])$ by $L_n(t)$. Then
 $L_n^b(t):=L_n(t\wedge\tau_n^b)$ is also a complex-valued
 martingale. It follows from the proof of Lemma \ref{lemma2.3} and Lemma
 \ref{lemma2.5} that
  \beqlb
  n(e^{-i\lambda/c_n}g_n(e^{i\lambda/c_n})-1)\rightarrow0
  \quad\mbox{and}\quad n^{\frac{1}{2}}(h_n(e^{i\lambda/c_n})-1)\rightarrow
  i\omega\gamma_0^{\frac{1}{2}}\lambda,
  \label{I_n}
  \eeqlb
as $n\rightarrow\infty$. Then we have for sufficiently large $n$,
 \beqlb
 \ln^{A_n(x_1,x_2,\lambda)}\ar=\ar i\lambda(m_n-1)x_1 +
 (c_nx_1+nx_2)\int_{\Gamma_n}(e^{i\lambda u}-1-i\lambda u)\mu_n(du)\nonumber\\
 \ar \ar \,+ (c_nx_1+nx_2)I_{1,n}(\lambda) +
 \int_0^{\infty}(e^{i\lambda u}-1-i\lambda u)\nu_n(du)\nonumber\\
 \ar \ar \, -\frac{1}{2}\big(h_n(e^{i\lambda/c_n})-1\big)^2 +
 I_{2,n}(\lambda),
 \label{taylor}
 \eeqlb
where
 \beqnn
 I_{1,n}(\lambda) \ar=\ar
 \sum_{j=2}^{\infty}(-1)^{j-1}\frac{[e^{i\lambda/c_n}g_n(e^{i\lambda/c_n})-1]^j}{j},\\
 I_{2,n}(\lambda)\ar=\ar\sum_{j=3}^{\infty}(-1)^{j-1}\frac{[h_n(e^{i\lambda/c_n})-1]^j}{j}.
  \eeqnn
Note that $n^2|I_{1,n}(\lambda)|\leq
\big|n(e^{i\lambda/c_n}g_n(e^{i\lambda/c_n})-1)\big|^2\rightarrow0$
and $n|I_{2,n}(\lambda)|\rightarrow0$. By {\rm (i)-(iii)}, Lemma
\ref{lemma2.5}, (\ref{I_n}) and (\ref{taylor}), it is not hard to
show that $n(A_n(x_1,x_2,\lambda)-1)\rightarrow A(x_1,x_2,\lambda)$
locally uniformly on $(x_1,x_2)\in\mathbb{R}\times\mathbb{R}_+$ for
fixed $\lambda$. As in  Ethier and Kurtz \cite[Problem 26,
P153]{EK86}, we obtain that
 \beqnn
\int^{t}_0e^{i\lambda
 Z^b_n(s)}n\big[A_n\big(Z^b_n(s),\phi^b_n(s),\lambda\big)-1\big]ds\rightarrow
 \int_0^t e^{i\lambda
 Z^b(s)}A(Z^b(s), \phi^b(s), \lambda)ds
 \eeqnn
in the topology of $C([0,\infty),\mathbb{C})$. Let
$L^b(t)=L(t\wedge\tau^b)$. Note that $[nt]/n\rightarrow t$ in
$C([0,\infty),\mathbb{R}_+)$. By (\ref{stopping times}),
\cite[Problem 13, P.151]{EK86} and \cite [Proposition 1.23,
p.293]{JS87}, we have
 \beqlb
 L_n^b(t)\overset{a.s.}{\longrightarrow}L^b(t)\quad \mbox{in}\
 D([0,\infty), \mathbb{C}),\ \quad \mbox{as}\ n\rightarrow\infty.
 \label{equation 5.14}
 \eeqlb
Then for almost all $t\geq0$,
$L^b_n(t)\overset{a.s.}{\longrightarrow}L^b(t)$ in $\mathbb{C}$. Fix
arbitrary $T>0$. For any $t\leq T$,
 $
\big|\int^{\tau_n^b(t)}_0e^{i\lambda Z^b_n(s)}\,Z_n^b(s) ds\big|\leq
bT,
 $
where the bound holds uniformly in $n$. Then for almost $t\leq T$,
$L^b_n(t)\overset{L_1}{\longrightarrow}L^b(t)$, as
$n\rightarrow\infty$. Since $L^b(t)$ is right continuous and bounded
for $t\leq T$, we have $L^b(t)$ is a martingale. Note that
$\tau^b\rightarrow\infty$ as $b\rightarrow\infty$, $L(t)$ is a local
martingale. \qed

It follows from (\ref{martingale problem}) and \cite [Theorem
2.42]{JS87} that $Z(\cdot)$ is a semimartingale and it admits the
canonical representation
 \beqlb
 Z(t)\ar=\ar Z(0) + Z^c(t)+\int_0^t\big(\beta_2+\beta_1\phi(s)+aZ(s)\big)\,ds +
 \int_0^t\int_0^\infty u\tilde{J}(ds,du),
 \label{semimartingale representation}
 \eeqlb
 where  $Z(0)=0$,
$Z^c(t)$ is a continuous local martingales with quadratic
covariation process $\int_0^t\varrho(s)ds$ with $
\varrho(s)=(2\sigma_1-a\gamma_0)\phi(s)+2\sigma_2+\omega(1-\omega)\gamma_0$,
and $J(dt,dz)$ is an integer-valued random measure on
$(0,\infty)\times\mathbb{R}_+$ with compensator
$\hat{J}(dt,du)=\phi(t)dt\mu(du)+dt\nu(du)$, where
$\tilde{J}(dt,dz)=J(dt,dz)-\hat{J}(dt,du)$.
 \blemma\label{lemma2.9} Suppose that the conditions of Theorem \ref{OU limit} are satisfied. Then the c\`{a}dl\`{a}g process
$Z(\cdot)$ is a weak solution of (\ref{OU}).
 \elemma
\proof Define the measure $\rho(du,d\zeta)=
\mu(du)\iota(d\zeta)+\nu(du)\delta_0(d\zeta)$, where $\iota(d\zeta)$
is the Lebesgue measure on $(0,\infty)$ and $\delta_0(d\zeta)$ is
the Dirac measure at $\zeta=0$. By Ikeda and Watanabe \cite[P.84 and
P.93]{IW89}, there exists a standard extension of $(\Omega,
\mathcal{F}, \mathcal{F}_t, P)$ supporting a one-dimensional
Brownian motion and a Poisson random measure $N(dt,du,d\zeta)$ on
$(0,\infty)\times\mathbb{R}^2_+$ with intensity $ds\rho(du,d\zeta)$
such that $dZ^c(t)=\sqrt{\varrho(t)}dB(t)$, and
 \beqlb
J((0,t]\times E)=\int_0^t\int_{\mathbb{R}^2_+}
1_{E}\big(\tilde{\theta}(s,u,\zeta)\big) N(ds,du,d\zeta),
 \eeqlb
for any $E\in\mathfrak{B}(\mathbb{R}_+)$, where
$\tilde{\theta}(s,u,\zeta)=u1_{[0,\phi(s)]}(\zeta)$. Set
$N_0(ds,du)=N(ds,du,\{0\})$ and set $N_1(ds,du,d\zeta)$
$=N(ds,du,d\zeta)|_{(0,\infty)\times\mathbb{R}_+\times(0,\infty)}$.
Then we see that $Z(\cdot)$ is a solution of (\ref{OU}). \qed

{\it Proof of Theorem \ref{OU limit}\;} By \cite[P.231]{IW89}, the
Lipschitz conditions of the equation (\ref{OU}) imply its pathwise
uniqueness of solutions. Thus Theorem \ref{OU limit} follows from
Lemma \ref{lemma2.7} and \ref{lemma2.9}.\qed

{\it Proof of Theorem \ref{asymptotic2.3}\;} By Theorem \ref{Theorem
1} $\big(Y_n(\cdot)\big/b_n, Y'_n(\cdot)\big/b_n\big)$ converges
weakly to $\big(Y(\cdot), Y'(\cdot)\big)$ on $D([0,\infty),
\mathbb{R}_+^2)$, and $\big(Y(\cdot), Y'(\cdot)\big)$ is
stochastically continuous. Note that
  \beqnn
  n(\check{m}-m_n)=\frac{Y_n(1)/b_n-n(m_n-1)\int_0^1Y_n(t)/b_n\,dt-Y'_n(1)/b_n}
  {\int_0^1Y_n(t)/b_n\,dt}.
  \eeqnn
Then we have (\ref{nature}) by the continuous mapping theorem. Since
$F$ is not bounded, the immigration process $Y'(\cdot)$ is neither a
compound Poisson process or a zero process. This implies that
$P(Y(t)=0\ \mbox{for all}\ t\in[0,1])=0$. \qed
 \blemma\label{lemma3.1}\;Define
 $u_n(k)=y_n(k)-m_ny_n(k-1)-\omega_n$. Then we have
  \beqlb
  \frac{1}{nc_n}\sum_{k=1}^nu_n^2(k)\overset{p}{\longrightarrow}0,\quad \mbox
  {as}\  n\rightarrow\infty.
  \label{limit0}
  \eeqlb
 \elemma
\proof It follows from (\ref{GWI-equationn}) and (\ref{regression
equation}) that
 \beqnn
 u_n(k)=\sum_{j=1}^{y_n(k-1)}(\xi_n(k,j)-m_n)+(\eta_n(k)-\omega_n).
 \eeqnn
Recall that $\hat{\xi}_n(k,j)$ and $\hat{\eta}_n(k)$ defined in the
proof of Lemma \ref{lemma2.6}. Note that
$\hat{\xi}_n(k,j)-\chi(\hat{\xi}_n(k,j))\geq0$ and
$\hat{\eta}_n(k)-\chi(\hat{\eta}_n(k))\geq0$. Then we have
  \beqlb
\frac{1}{nc_n}\sum_{k=1}^nu_n^2(k)\ar\leq\ar6(I_{1,n}^2 + I_{2,n})+
\frac{6c_n}{n}\Big[n^{\frac{3}{2}}\int_{\Gamma}(\chi(u)-u)\mu_n(du)\Big]^2
\int_0^1\Big(\frac{Y_n(s)}{n}\Big)^2ds \nonumber\\
 \ar \ar
+\, 6c_n\Big(\int_0^\infty(\chi(u)-u)\nu_n(du)\Big)^2 +
6\Big[\sqrt{\frac{c_n}{n}}\,\sum_{k=1}^n\big(\hat{\eta}_n(k)-\chi(\hat{\eta}_n(k))\big)\Big]^2
  \nonumber\\
  \ar \ar
+\,\frac{6c_n}{n}\sum_{k=1}^n(\chi(\hat{\eta}_n(k))-\bfE[\chi(\hat{\eta}_n(k))])^2,
   \label{W_n(k)0}
    \eeqlb
where
 \beqnn
 I_{1,n}\ar=\ar\sqrt{\frac{c_n}{n}}\,\sum_{k=1}^n\sum_{j=1}^{y_n(k-1)}
 \big[\hat{\xi}_n(k,j)-\chi(\hat{\xi}_n(k,j))\big],\\
 I_{2,n}\ar=\ar
 \frac{c_n}{n}
 \sum_{k=1}^n\Big[\sum_{j=1}^{y_n(k-1)}\big(\chi(\hat{\xi}_n(k,j))-\bfE[\chi(\hat{\xi}_n(k,j))]\big)\Big]^2.
 \eeqnn
We obtain that $\bfE[I_{1,n}]\leq\sqrt{\frac{c_n}{n}}\Phi(1)K$ and
$\bfE[I_{2,n}]\leq\frac{c_n}{n}\Phi(1)\big[n^2\,\bfvar\,\chi\big({\hat{\xi}_n(1,1)}\big)\big]$
as in the calculations in (\ref{Z_ninequality2}). Condition (C)
implies that $I_{i,n}\overset{p}{\longrightarrow}0$ as
$n\rightarrow\infty$ for $i=1,2$. From (C), Remark \ref{remark 2.1},
and {\rm (a.1)} in the proof of Lemma \ref{lemma2.3}, the third term
in (\ref{W_n(k)0}) converges in probability to $0$ as
$n\rightarrow\infty$. As in the above proof, we also have that the
last three terms converge in probability to $0$. Thus
(\ref{W_n(k)0}) implies (\ref{limit0}).\qed

 {\it Proof of Theorem \ref{asymptotic}\;} First consider the
 equation (\ref{OU}). We obtain as in the calculation in
 (\ref{Z_ninequality2}) and (\ref{Aldous}) that for $0\leq s\leq t$,
  \beqlb
  \bfE[|Z(t)-Z(s)|]\ar\leq\ar|\beta_2|(t-s)
  +|\beta_1|\int_s^t\phi(u)du+|a|\int_s^t\hat{\Phi}(u)e^{|a|u}du\nonumber\\
   \ar \ar
   +\, 2\int_s^t\phi(u)du\int_0^\infty(u\wedge u^2)\mu(du)
   +2(t-s)\int_0^\infty(u\wedge u^2)\nu(du)
   \nonumber\\
   \ar \ar
   +\,\Big(\int_s^t\phi(u)du\int_0^\infty(u\wedge
   u^2)\mu(du)\Big)^{\frac{1}{2}} +
   \Big(\int_s^t\varrho(u)du\Big)^{\frac{1}{2}},
  \eeqlb
where $\hat{\Phi}(\cdot)$ is some non-decreasing continuous
function. Then $Z(\cdot)$ is stochastically continuous. Let
$D(Z):=\{t\geq0: \bfP\{Z(t)=Z(t-)\}=1\}$ and thus $D(Z)=(0,\infty)$.
From (\ref{form}) we obtain
 \beqlb
 \frac{n^2}{c_n}(\hat{m}_n-m_n)=\frac{\frac{1}{nc_n}\sum_{k=1}^ny_n(k-1)u_n(k)}
 {\frac{1}{n^3}\sum_{k=1}^ny_n^2(k-1)}=\frac{D(n)}{Q(n)}.
 \label{form2}
 \eeqlb
Rewrite $D(n)$ as
$D(n)=D_1(n)+\frac{c_n}{n}\sum_{j=2}^3D_j(n)-D_4(n)$, where
 \beqlb
 D_1(n)\ar=\ar \frac{1}{nc_n}\sum_{k=1}^n\bfE[y_n(k-1)]u_n(k),
\quad D_2(n)=\frac{n(1-m_n^2)}{2m_n}\int_0^1Z_n(s)ds,\nonumber\\
 D_3(n)\ar=\ar\frac{1}{2m_n}Z^2_n(1),\qquad\quad
 D_4(n)=\frac{1}{2nc_nm_n}\sum_{k=1}^nu_n^2(k).
 \eeqlb
Let $M_n(t)=\sum_{k=1}^{[nt]}u_n(k)/c_n$ for $t\geq0$. The
functional $\Psi_n: D([0,\infty),\mathbb{R})\mapsto \mathbb{R}$ is
defined by
 \beqlb
 \Psi_n(x)=\omega_n\int_0^1(x(1)-x(t))m_n^{[nt]-1}dt-\frac{x(1)}{nc_nm_m}.
 \eeqlb
Then $D_1(n)$ can be rewritten as
   \beqnn
D_1(n)=\frac{1}{nc_n}\sum_{j=1}^{n-1}m_n^{j-1}\sum_{k=j+1}^nu_n(k)=\Psi_n(M_n).
   \eeqnn
If $x_n\rightarrow x$ in the topology of $D([0,\infty),\mathbb{R})$,
it is easy to see that $|\Psi_n(x_n)-\Psi(x)|\rightarrow0$, where
  \beqnn
\Psi(x)=\omega\int_0^1(x(1)-x(t))e^{at}dt.
  \eeqnn
Note that $M_n(t)=Z_n(t)-\int_0^{[nt]/n}n(m_n-1)Z_n(s)ds$ by
(\ref{GW difference equation}). It follows from Theorem \ref{OU
limit} that $M_n(t)$ converges weakly to ${M}(t):=Z(t)-\int_0^t
aZ(s)ds$ on $D([0,\infty),\mathbb{R})$. By Remark \ref{remark 2.1}
$Y_n(\cdot)/n$ converges weakly to $\phi(\cdot)$ on
$D([0,\infty),\mathbb{R}_+)$, and $\phi(\cdot)$ is a deterministic
continuous function. Thus $\big(Y_n(\cdot)/n, Z_n(\cdot),
M_n(\cdot)\big)$ converges weakly to $\big(\phi(\cdot), Z(\cdot),
M(\cdot)\big)$ on $D([0,\infty),\mathbb{R}_+\times\mathbb{R}^2)$. By
\cite[Theorem 7.8, P.131 and Problem 26, P.153]{EK86} and the
continuous mapping theorem, we have
$D_1(n)\overset{d}{\rightarrow}\Psi(M)=\int_0^1\phi(s)dM(s)$,
$D_2(n)\overset{d}{\rightarrow} -\int_0^1aZ(s)ds$,
$D_3(n)\overset{d}{\rightarrow} \frac{1}{2}Z^2(1)$ and
$Q_n\overset{p}{\rightarrow}\int_0^1\phi(s)ds$, as
$n\rightarrow\infty$. Then it follows from (C) and Lemma
\ref{lemma3.1} that
$\frac{c_n}{n}\sum_{j=2}^3D_j(n)-D_4(n)\overset{p}{\rightarrow}0$ .
Hence we obtain (\ref{asymptotic1}). In a similar way, we also have
(\ref{asymptotic1+}).\qed

Recall that $u_n(k)=y_n(k)-m_ny_n(k-1)-\omega_n$. Let
$v_n(k)=u_n^2(k)-\pi_ny_n(k-1)-r_n$ and  let
$V_n(t)=\sum_{k=1}^{[nt]}v_n(k)$. By (\ref{variance equation}),
$V_n(\cdot)$ is a martingale.

{\it Proof of Theorem \ref{asymptotic variance}\;} Under the
conditions of Corollary \ref{example1}, $Z_n(\cdot)$ is defined by
(\ref{Z_n}) with $c_n=\sqrt{n}$, and then
$M_n(t)=\sum_{k=1}^{[nt]}u_n(k)/\sqrt{n}$. By Corollary
\ref{example1} and the proof of Theorem \ref{asymptotic}, we have
that $(Z_n(\cdot), M_n(\cdot))$ converges weakly to $\big(Z(\cdot),
M(\cdot)\big)$ on $D([0,\infty),\mathbb{R}^2)$, where $Z(\cdot)$ is
given by (\ref {OU_square equation}) and ${M}(t)=Z(t)-\int_0^t
aZ(s)ds$. Note that $M_n(\cdot)$ is a square integrable martingale
and
$\bfE[M_n^2(t)]=(\pi_n\sum_{k=1}^{[nt]}\bfE[y_n(k-1)]+[nt]r_n)/n$.
Then for $t\geq0$ and sufficiently large $n$,
$\bfE[M_n^2(t)]\leq\tilde{\mu}(\mathbb{R}_+)
\int_0^t\phi(s)ds+\tilde{\nu}(\mathbb{R}_+)t+1$. Thus by Kurtz and
Protter \cite[Theorem 2.7]{KP91},
 \beqlb
 \bigg(Z_n(t), M_n(t),
\int_0^tZ_n(s-)dM_n(s)\bigg)\rightarrow \bigg(Z(t), M(t),
\int_0^tZ(s-)dM(s)\bigg)
 \label{integral convergence}
 \eeqlb
in distribution on $D([0,\infty),\mathbb{R}^3)$. On the other hand,
let $\hat{V}_n(t):=V_n(t)/n$.
  \beqnn
\hat{V}_n(t)&=&Z_n^2(t)+n(1-m_n^2)\int_0^{[nt]/n}Z_n^2(s)ds
-2m_n\int_0^tZ_n(s-)dM_n(s)\nonumber\\
 & &\;-n\pi_n\int_0^{[nt]/n}Y_n(s)/n\,ds-[nt]r_n/n.
   \eeqnn
Still note that $Y_n(\cdot)/n$ converges weakly to $\phi(\cdot)$ on
$D([0,\infty),\mathbb{R}_+)$, and $\phi(\cdot)$ is a deterministic
continuous function. By (\ref{integral convergence}) and the
continuous mapping theorem, $(Z_n(\cdot),\hat{V}_n(\cdot))$
converges weakly to $(Z(\cdot), J(\cdot))$ on
$D([0,\infty),\mathbb{R}^2)$, where
$J(t)=Z^2(t)-2a\int_0^tZ^2(s)ds-2\int_0^tZ(s)dM(s)
-\tilde{\mu}(\mathbb{R}_+)\int_0^t\phi(s)ds-\tilde{\nu}(\mathbb{R}_+)t$.
By It\^{o}'s formula, $J(t)$ has also the form (\ref{J(t)}).
$\hat{V}_n(t)$ is also a finite variation process. Denote its finite
variation by $\int_0^t|dV_n(s)|$. Then for $t\geq0$ and sufficiently
large $n$,
$\bfE\big[\int_0^t|d\hat{V}_n(s)|\big]\leq2\tilde{\mu}(\mathbb{R}_+)
\int_0^t\phi(s)ds+2\tilde{\nu}(\mathbb{R}_+)t+1$. By \cite{KP91}
again, $\int_0^tY_n(s-)/n\,d\hat{V}_n(s)$ converges weakly to
$\int_0^t\phi(s)dJ(s)$ on $D([0,\infty),\mathbb{R})$. We see that
 \beqlb
 n(\hat{\pi}_n-\pi_n)=\frac{\int_0^1Y_n(s-)/n\,d\hat{V}_n(s)-\hat{V}_n(1)\int_0^1Y_n(s)/n\,ds}
 {\int_0^1(Y _n(s)/n-\int_0^1Y_n(s)/n\,ds)^2ds},
  \label{calculation of pi}
 \eeqlb
and
$\hat{r}_n-r_n=\hat{V}_n(1)-(\hat{\pi}_n-\pi_n)\int_0^1Y_n(s)ds$.
Note that $J(t)$ and $\int_0^t\phi(s)dJ(s)$ are stochastically
continuous. By the continuous mapping theorem, we have (\ref
{asymptotic2+}). We write
  \beqnn
u_n(k)-\hat{u}_n(k)&=&-[(\hat{m}_n-m_n)y_n(k-1)]^2-(\hat{\omega}_n-\omega_n)^2
+2(\hat{m}_n-m_n)y_n(k-1)u_n(k)\\
 &&+2(\hat{\omega}_n-\omega_n)u_n(k)-2(\hat{m}_n-m_n)(\hat{\omega}_n-\omega_n)y_n(k-1).
  \eeqnn
As in the proof of (\ref{calculation of pi}), also by Theorem
\ref{asymptotic}, we have (\ref {asymptotic2+}) holds for
$\tilde{\pi}_n$ and $\tilde{r}_n$. \qed

By the proof of Theorem \ref{asymptotic variance}, we see that
$\hat{V}_n(t):=V_n(t)/n$ converges weakly to $J(\cdot)$ on
$D([0,\infty),\mathbb{R})$, where $J(t)$ is defined by (\ref{J(t)}).
However when we turn to the case of Example \ref{example of
diffusion}, $J(t)$ is degenerate to $0$. Then in this case, we need
the following lemma.

 \blemma\; Let $\bar{V}_n(t)=V_n(t)/\sqrt{n}$. Under
 the conditions of Theorem \ref{diffusion variance theorem}, $\bar{V}_n(\cdot)$ converges in
distribution on $D([0,\infty), \mathbb{R})$ to the process
$V(\cdot)$, which is defined by
$V(t)=\int_0^t\sqrt{\varrho_2(s)}dW(s)$, where
$\varrho_2(t)=2\pi^2\phi^2(t)+(a_4+4\pi r)\phi(t)+(b_4-r^2)$ and
$W(\cdot)$ is a one-dimensional Brownian motion.
 \elemma
\proof Under the above conditions, we see that $Y_n(\cdot)/n$
converges weakly to $\phi(\cdot)$ on $D([0,\infty),\mathbb{R}_+)$ by
Remark \ref{remark 2.1}  {\rm(ii)}. Then for any $t\geq0$,
  \beqnn
  \frac{1}{n}\sum_{k=1}^{[nt]}\bfE[v_n^2(k)|\mathcal{F}_{k-1}^n]
  =
 \frac{1}{n}\sum_{k=1}^{[nt]}\Big(
  (a_{4,n}+4\pi_nr_n-3\pi_n^2)y_n(k-1)+2\pi_n^2y_n(k-1)^2+b_{4,n}-r_n^2\Big),
  \eeqnn
which converges in probability to $\int_0^t\varrho_2(s)ds$. Now by
the martingale central limit theorem, it suffices to prove that, for
any $\varepsilon>0$ and $t\geq0$,
 \beqlb
 \frac{1}{n}\sum_{k=1}^{[nt]}\bfE[v_n^2(k)1_{\{|v_n(k)|>\sqrt{n}\varepsilon\}}
 |\mathcal{F}_{k-1}^n]\overset{p}{\rightarrow}0, \quad\mbox{as}\
 n\rightarrow\infty.
 \label{lindeberg}
 \eeqlb
We have that $v_n(k)=A_{n,k}+B_{n,k}+C_{n,k}+D_{n,k}$, where
 \beqnn
\left.\begin{array}{ll}\displaystyle
A_{n,k}=\sum_{i=1}^{y_n(k-1)}[(\xi_n(k,i)-m_n)^2-\pi_n], \quad
B_{n,k}=2\sum_{i=1}^{y_n(k-1)}(\xi_n(k,i)-m_n)(\eta_n(k)-\omega_n), \\
\displaystyle C_{n,k}=(\eta_n(k)-\omega_n)^2-r_n,\quad
D_{n,k}=2\sum_{i<j}^{y_n(k-1)}(\xi_n(k,i)-m_n)(\xi_n(k,j)-m_n).
  \end{array}
  \right.
  \label{eq2.6}
  \eeqnn
Note that for any pair of random variables $\bar{X}$ and $\bar{Y}$,
$\bfE\big[(\bar{X}+\bar{Y})^21_{\{|\bar{X}+\bar{Y}|>\varepsilon\}}\big]\leq
4\Big(\bfE\big[\bar{X}^21_{\{|\bar{X}|>\varepsilon/2\}}\big]
+\bfE\big[\bar{Y}^21_{\{|\bar{Y}|>\varepsilon/2\}}\big]$\Big). Thus
it suffices to show that (\ref{lindeberg}) with $v_n(k)$ replaced by
$A_{n,k}$, $B_{n,k}$, $C_{n,k}$, and $D_{n,k}$. Let
$\xi'_n(k,i)=(\xi_n(k,i)-m_n)^2-\pi_n$. As in the proof of
\cite[Theorem 2.2]{I05}, we obtain
  \beqnn
 &&\frac{1}{n}\sum_{k=1}^{[nt]}\bfE[A_n^2(k)1_{\{|A_n(k)|>\sqrt{n}\varepsilon\}}
 |\mathcal{F}_{k-1}^n] \\
 &\leq&
 n\bfE\big[(\xi'_{n}(1,1))^21_{\{|\xi'_n(1,1)|>\sqrt{n}\varepsilon/2\}}\big]
  \sum_{k=1}^{[nt]}y_n(k-1)/n^2
  +4n(a_{4,n}-\pi_n^2)^2\varepsilon^{-2}
  \sum_{k=1}^{[nt]}y_n^2(k-1)/n^3\\
  &&+\sqrt{2}n(a_{4,n}-\pi_n^2)^{\frac{3}{2}}\varepsilon^{-1}
  \sum_{k=1}^{[nt]}y_n^{\frac{3}{2}}(k-1)/n^{\frac{5}{2}}.
  \eeqnn
For large enough $n$, $\pi_n\leq\sqrt{n}\varepsilon/2$, and
 \beqnn
\bfE\big[(\xi'_{n}(1,1))^21_{\{|\xi'_n(1,1)|>\sqrt{n}\varepsilon\}}\big]
&\leq&\bfE\big[(\xi_n(1,1)-m_n)^4{\bf1}_{\{(\xi_n(1,1)-m_n)^2>\sqrt{n}\,\varepsilon/2
 \}}\big]\\
  &&+2(\pi_na_{4,n}+2\pi_n^3)/\sqrt{n}\varepsilon.
  \eeqnn
Then, by conditions (c.1,2), (\ref {lindeberg}) holds with $v_n(k)$
replaced by $A_{n,k}$. Also by condition (c.3), (\ref {lindeberg})
holds for $C_{n,k}$. Let $\bar{\xi}_n(k,i)=\xi_n(k,i)-m_n$. For
$D_{n,k}$, we note that
  \beqnn
  D_{n,k}^2/4&=&\pi_n\sum_{j=2}^{y_n(k-1)}(j-1)\xi'_n(k,j)
  +\pi_n\sum_{i=1}^{y_n(k-1)-1}(y_n(k-1)-i)\xi'_n(k,i)\\
  &&+\,\sum_{i<j}^{y_n(k-1)}\xi'_n(k,i)\xi'_n(k,j)+2\sum_{l<i<j}^{y_n(k-1)}
  (\bar{\xi}_n(k,l))^2\bar{\xi}_n(k,i)\bar{\xi}_n(k,j)\\
  &&+2\sum_{l<i<j}^{y_n(k-1)}
  \bar{\xi}_n(k,l)(\bar{\xi}_n(k,i))^2\bar{\xi}_n(k,j)+
  2\sum_{l<i<j}^{y_n(k-1)}
  \bar{\xi}_n(k,l)\bar{\xi}_n(k,i))(\bar{\xi}_n(k,j))^2\\
  &&+6\sum_{l<i<j<p}^{y_n(k-1)}
  \bar{\xi}_n(k,l)\bar{\xi}_n(k,i)\bar{\xi}_n(k,j)\bar{\xi}_n(k,p)
  +y_n(k-1)(y_n(k-1)-1)\pi_n^2/2.
  \eeqnn
Then it follows from the above equality that
  \beqlb
  \bfE[D_{n,k}^4|\mathcal{F}_{k-1}^n]\leq16a_{4,n}^2y_n^2(k-1)
  +416a_{4,n}\pi_n^2y_n^3(k-1)
   +772\pi_n^4y^4_n(k-1).
   \label{B_n fourth moment}
  \eeqlb
Thus, for any $t\geq0$,
$\frac{1}{n}\sum_{k=1}^{[nt]}\bfE[D_n^2(k)1_{\{|D_n(k)|>\sqrt{n}\varepsilon\}}
 |\mathcal{F}_{k-1}^n]\leq\frac{1}{n^2\varepsilon^2}
 \sum_{k=1}^{[nt]}\bfE[D_n^4(k)|\mathcal{F}_{k-1}^n]$, which
 converges in probability to $0$ by (\ref{B_n fourth moment}). In a
 similar way, we can also prove that (\ref{lindeberg}) holds with $v_n(k)$ replaced by
 $B_{n,k}$.\qed

{\it Proof of Theorem \ref{diffusion variance theorem}\;} It is not
hard to see that for any $t\geq0$,
  \beqnn
  \frac{1}{n^3}\sum_{k=1}^{[nt]}E[y^2_n(k-1)]&=&\frac{\pi_n+2\omega_nm_n}{m^2_n}
  \int_0^{\frac{[nt]}{n}}m^{2[ns]}_n\int_0^{\frac{[ns]}{n}}m^{-2[nu]}_n\int_0^{\frac{[nu]}{n}}m^{[n\zeta]}d\zeta\,
  \,du\,ds\\
   &&+\,\frac{r_n+\omega_n^2}{n}\int_0^{\frac{[nt]}{n}}\int_0^{\frac{[ns]}{n}}m^{2[nu]}du\,ds,
  \eeqnn
which converges to $\int_0^t\phi^2(s)ds$ as $n\rightarrow\infty$.
Also by (\ref{expectation limit}) and the proof of Lemma \ref
{lemma2.7}, we see that $\bar{V}_n(\cdot)$ is a square integrable
martingale and $\sup_n\bfE[\bar{V}_n^2(t)]<\infty$ for $t\geq0$. By
\cite[Theorem 2.7]{KP91},
 $
 \big(Y_n(t)/n, \bar{V}_n(t),
\int_0^tY_n(s-)/n\,d \bar{V}_n(s)\big)\rightarrow \big(\phi(t),
V(t), \int_0^t\phi(s)dV(s)\big)
 $
in distribution on $D([0,\infty),\mathbb{R}_+\times\mathbb{R}^2)$.
As in the proof of Theorem \ref{asymptotic variance}, we have
(\ref{asymptotic2+2}).\qed

\bigskip
\noindent\textbf{Acknowledgement.} I would like to thank my
supervisor Professor Zenghu Li for his encouragement and many
fruitful discussions. This work was supported by NSFC
(No.\,10871103).

\end{document}